%% file: enumii.tex
\documentclass[twoside]{article}
\usepackage{amssymb, floatflt}
\input enum.mac
\input enumi.aux

\renewcommand{\copyrightheading}[1]{}

\runninghead{Enumerating the Prime Alternating Knots, Part II }
{Enumerating the Prime Alternating Knots, Part II }
\begin{document}
\copyrightheading{}		    

\fpage{1}
\centerline{\bf Enumerating the Prime Alternating Knots, Part II}
\centerline{\footnotesize Stuart Rankin, John Schermann, Ortho Smith\fnm{$\dagger$}}
\fnt{$\dagger$}{Partially supported by the Applied Mathematics Department at the University of
Western Ontario, Wolfgang Struss, Pei Yu, Turab Lookman, Ghia Flint and Dean Harrison.}
\centerline{\footnotesize Department of Mathematics, University of Western Ontario}
\centerline{\footnotesize srankin@uwo.ca, johns@csd.uwo.ca, dsmith6@uwo.ca}
\vskip20pt

\input introii.tex\par
\input master.tex\par
\input implem1.tex\par
\input implem2.tex\par
\input implem3.tex\par

\section{References}

[1] S. A. Rankin, J. Schermann, O. Smith, Enumerating the prime alternating
 knots, Part I, preprint.

\end{document}

%% file: introii.tex
\section{Introduction}
 In Part I, we provided a solution for part of the enumeration problem,
 in that we described an inductive scheme which used a total of
 four operators to generate all prime alternating knots of a given
 minimal crossing size from those of crossing size one less, and we provided
 a proof there that the procedure does produce them all. Our primary
 objective in Part II is to present an efficient implementation of this
 procedure.

 One very important aspect of our implementation is the notion of the
 master array of an alternating knot. The master array of a prime
 alternating knot is an integer array with the property that
 each regular projection (plane configuration) at minimal crossing size
 for the knot can be constructed from the data in the array, and moreover,
 two knots are flype equivalent if and only if their master arrays are
 identical. Furthermore, the master array for a prime alternating knot can
 be calculated from the data describing any minimal crossing configuration
 for the knot, and we have a procedure for selecting from the master
 array a preferred configuration, which we suggest could be used as the
 so-called ideal configuration of the knot. A very efficient way to store
 the knot would then be to use the Dowker-Thistlethwaite code for this
 ideal configuration.
  
 An early implementation of this algorithm was used to enumerate the prime
 alternating knots up to and including those of 19 crossings. It took
 approximately $2.3$ hours of time on a five node beowulf cluster to
 produce the $1,769,979$ prime alternating knots of 17 crossings. We
 then went on to produce the prime alternating knots at 18 and 19 crossings
 using a 48 node beowulf cluster. The cluster was shared with other users
 and so an accurate estimate of the running time is not available, but
 the generation of the $8,400,285$ knots at 18 crossings was completed in
 17 hours, and the
 generation of the $40,619,385$ prime alternating knots at 19 crossings
 took approximately 72 hours. In Part II, we describe improvements to our
 original implementation that should allow us to produce the prime alternating
 knots at 19 crossings from those at 18 crossings in about 10 hours on a
 current Pentium III personal computer equipped with 256 megabytes of
 main memory.

Since we shall only be working with minimal crossing configurations of
alternating knots, it should be understood from now on that when we
say configuration, we mean minimal crossing configuration.
 

%% file: master.tex
\def\lmaxs{\ell_{\hbox{max}}}
\def\lmax#1{\lmaxs(#1)}
\def\Lmax{L_{\hbox{max}}}

\section{Master Group Codes and the Master Array}

 In this section, we introduce the master array for a prime
 alternating knot. The master array is constructed from any
 group code for any configuration of the knot. It describes the
 basic structure of the knot and identifies all possible flype
 moves for the knot. Consequently, it is possible to construct
 every configuration of the knot from its master array. Moreover,
 two knot configurations are flype equivalent if and only if the
 master arrays that are constructed from the two configurations
 are identical.

 As a first step toward the construction of the master array for a
 knot, we introduce the notion of a master group code. A
 master group code is constructed from any configuration for
 the knot. It is a cyclic code which contains the complete orbit
 structure of the knot, and therefore any configuration of the knot
 can be constructed from any master group code for the knot. There
 are, however, many master group codes for a given knot, and this
 would lead to extensive computation if one were to try to compare two
 configurations for flype equivalency by comparing  master group
 codes for the two configurations. The master array is actually
 a master group code which satisfies certain conditions (which are such
 that there is exactly one master group code that satisfies the conditions),
 together with a pointer to the position from which the traversal of the
 code is to begin, and a direction of traversal.

 We begin with the definition of a master group code for
 an alternating knot. Following that, we present an algorithm for the
 construction of all master group codes of a given knot. The section
 concludes with the definition and construction of the master array for
 a knot.

\subsection{The definition of master group codes}

 We have observed that the group orbits of a prime alternating knot
 can be determined from any  configuration of the knot.
 Furthermore, as established in Theorem \ref{full group} of [1],
 once the orbits of all groups are known, then it is a simple matter
 to convert the configuration into a full group configuration.
 More generally, once the groups and their orbits are identified, then
 every  configuration of the knot can be
 constructed by flyping the various crossings
 of the different groups, wherein each crossing can be placed at any
 position in the orbit of the group of the crossing. It turns out that
 the decision to deal only with full group configurations provided us 
 with a significant theoretical and computational advantage over the
 traditional methods in the study of the prime alternating knots.

 These observations lead naturally to the notion of what we call a
 master group code for the knot. In effect, a master group code for a knot
 describes the orbit structure of each group, the size of each group,
 and the sign (positive or negative) of each group.

 We shall need the following result, which is an immediate consequence
 of the non-interference of orbits (see Theorem \ref{non-interference
 of orbits} of [1]).

\begin{theorem}\label{non interference revisited}
 Let $C(K)$ be a full group configuration of a prime alternating knot
 $K$, and let $G_1$, $G_2$, $G_3$ and $G_4$ be distinct groups in $C(K)$.
 Further suppose that $G_3$ and $G_4$ are connected by an arc $e$ and
 that $G_1$ can group flype to $e$ and some other arc, and $G_2$ can
 group flype to $e$ and some other arc, so that $G_3$ and $G_4$ are in
 adjacent min-tangles $T_1$ and $T_2$, respectively, in the orbit of
 $G_1$, and they are in adjacent min-tangles $S_1$ and $S_2$, respectively,
 in the orbit of $G_2$. Then
 \begin{alphlist}
  \item $G_1$ is contained in either $S_1$ or $S_2$, and $G_2$ is contained
  in either $T_1$ or $T_2$;
  \item $G_1$ is in $S_1$ if and only if $G_2$ is in $T_2$.
  \item $G_1$ is in $S_1$ if and only if for any sequence of group
  flypes applied to $C(K)$ which results in both $G_1$ and $G_2$ being
  group flyped to lie on $e$, the knot traversal out from the group
  $G_3$ in the direction of $e$ encounters $G_1$ before $G_2$.
 \end{alphlist}
\end{theorem}

\begin{proof}
 $G_2$ belongs to some min-tangle of the orbit of $G_1$, and by
 Theorem \ref{non-interference of orbits} of [1], $G_2$ can't flype to
 any arcs that are not incident to or contained within this min-tangle.
 Since $G_2$ can flype onto arc $e$, which is incident to the two
 min-tangles $T_1$ and $T_2$ of the orbit of $G_1$, it follows that
 $G_2$ is contained either in $T_1$ or else $T_2$. Similarly, $G_1$
 belongs either to $S_1$ or to $S_2$.

 Suppose that $G_1$ is in $S_1$. Any group flype
 move applied to $G_2$ must leave it still in $T_1$. Since there is a
 group flype which will move $G_2$ onto $e$, let us apply it. Let
 $e'$ denote the arc connecting $G_3$ to $G_1$ and $e''$ denote the
 arc connecting $G_1$ to $G_4$ (in effect, $e$ has been replaced by
 $e'$ and $e''$). Now $e''$ is the arc connecting the min-tangle
 $T_1$ (or rather, the result of applying the group flype to $G_2$ in
 $T_1$) to the min-tangle $T_2$, and so we may group flype $G_1$ to
 $e''$ and the other arc joining $T_1$ to $T_2$. Now the orbits of
 $G_1$ and $G_2$ are not changed by the group flyping of $G_1$ and
 $G_2$, in the sense that the groups that make up a given min-tangle $T$
 of a given group are still the groups that make up the min-tangle
 after any group flype whatsoever--they may have changed location,
 but they have not left the min-tangle $T$. In this last configuration,
 we have $G_3$ connected to $G_1$, which in turn is connected to $G_2$.
 Thus $G_3$ is in one min-tangle of the orbit of $G_1$ and $G_2$ is in
 a different min-tangle of the orbit of $G_2$. Since $G_3$ is in
 $T_1$, we see that $G_2$ is not in $T_1$. Then, since we know that
 $G_2$ is in either $T_1$ or else $T_2$, we may conclude that
 $G_2$ is in $T_2$, as required. The symmetry in this argument allows
 us to conclude similarly that $G_2$ in $T_2$ implies that $G_1$ is in
 $S_1$.

 For the last part, we again suppose that $G_1$ is in $S_1$. Apply any
 sequence of group flypes which results in both $G_1$ and $G_2$ being
 group flyped to the arc $e$. Then as we traverse the knot from $G_3$
 in the direction of $e$, we remain in the min-tangle $S_1'$ of the
 orbit of $G_2$ which is the image of $S_1$ under this seqence of group
 flypes. Since $G_1$ lies on $e$, we will encounter $G_1$ before leaving
 $S_1'$, hence will encounter $G_1$ before encountering $G_2$. Then
 $G_1$ and $G_4$ lie in different min-tangles of the orbit of $G_2$.
 Since $G_4$ is in $S_2$, this implies that $G_1$ is in $S_1$, as
 required.
\end{proof}

The gist of Theorem \ref{non interference revisited} is that whenever two
or more groups have flype positions with one edge $e$ in common, and the
groups are group flyped to the position in their respective orbits at
which $e$ appears, then the order in which they appear on that strand
in any traversal of the knot is independent of the order in which the
group flyping took place.

\begin{definition}
 A {\it master group code} for a prime alternating knot $K$ is constructed
 from a full group  configuration $C(K)$ of $K$ according
 to the following procedure. Let $C$ be any group code for the
 configuration. Recall that in the construction (see Definition \ref{group
 code def} of [1]) of the group code, each group $G$ was
 assigned a label of the
 form $\pm m_n$, where $m$ denotes the size of the group and $n$ was an
 index used to distinguish between different groups of the same size.
 Each pair of consecutive group labels in the group code (recall, this
 is a cyclic code, so the last group label of the code is followed by the
 first group label of the code) represents the arc that joins these two
 groups. For each group $\pm m_n$, list the pairs of arcs
 which separate consecutive min-tangles in the orbit of the group
 (excluding the arcs which connect the group itself to the min-tangles
 on either side of the group). If the orbit is trivial, there will be only
 one min-tangle in the orbit and so the list will be empty, in which
 case the two occurrences of the labels $m_n$ are replaced by $m_n^0$, keeping
 whatever sign the label had originally. In the case when the orbit is
 non-trivial, there will be one or more pairs of arcs in the list. In this case,
 let $k$ denote the number of pairs in the list. Form the labels $m_n^i$, $i
 =0,\ldots, k$. Choose one of these $k+1$ labels and replace both occurrences
 of $m_n$ in the group code by the selected label $m_n^i$, keeping the
 sign of the label being replaced. Then in an arbitrary fashion, assign
 the remaining $k$ labels to the $k$ pairs of arcs in the list. For each
 arc of a pair, place the label assigned to the pair between the two
 consecutive groups in the group code which identify the arc. If the group
 is a negative group, both arcs should receive a minus sign as part of
 the assigned label. In the event that
 more than one label is inserted between two consecutive entries of the
 group code, the labels must be arranged in the uniquely determined order
 (see Theorem \ref{non interference revisited}) in which the groups would
 appear if all were group flyped to their respective orbit positions which
 involved the arc in question.

 The cyclic arrangement of labels that results is a {\it master group code}
 for the knot.
\end{definition}

 In practice, the cyclic arrangement is written as a sequence, with the
 understanding that the initial entry follows that last entry of the
 sequence.
 
 Since two full group configurations for the same knot only differ
 in the position in its orbit that each full group is placed, it
 follows that all full group configurations for a knot can be
 constructed from the data stored in any master group code for the knot
 (and, of course, all split group configurations can be so obtained as
 well).

\subsection{An Algorithm for the Construction of a Master Group Code}

It should be clear now that any two full group 
configurations will result in two master group codes which are essentially
the same, except for possibly the choice of labelling. There is, however,
the possibility that the two group codes that were used to construct the two
master group codes for the same knot arose from traversals of the knot in
opposite directions. Thus if we were to try
to compare two master group codes to determine whether or not they describe
the same knot, we not only need to rotate one code to compare to the other,
but we must try to decide whether or not we would reverse one before doing
the rotations for attempted pattern matching. For much of our work, any
master group code would suffice, but it will be necessary for us to compare
two configurations for flype equivalence, and that leads to the problems
described above. For this reason, it is important to assign to a given knot
a representative of the set of master group codes, together with a start
point and direction of traversal with the property that two configurations
are flype equivalent if and only if the two designated master group codes,
complete with start point and direction of traversal, are identical.

Given a master group code, it is a simple matter to indicate a start point
and direction of traversal of the code: we simply open the cyclic code at
the desired point and write the entries in a sequence to describe the
order of traversal. That is, we form a linear array. The problem that
remains is how to decide on an appropriate scheme for selecting from
all of the arrays that can be formed from all of the different master
group codes for the knot the unique one with the desired attributes.
We will leave this problem for the last subsection of the current section,
and concentrate now on an algorithm for the construction of the master
group codes for a knot.

Actually, since it does not take much more work to accomplish this, we
shall present an algorithm which takes
as its input the group code for a  knot configuration
(either split or full group) of a prime alternating knot, identifies
the orbit of each group, identifies all split groups and recombines
the split groups so as to obtain a full group configuration, all the
while constructing a master group code for the knot as well.

We start with an array that initially consists of a group code for a
 configuration of the prime alternating knot. In a
sense, the group code serves as a skeleton for the
master group code that is to be constructed. The first step is to make
a list of all the groups in the knot configuration.
Then for each group in the list, we identify the min-tangles of the orbit
of the group.
In each case, if the group is a split group, then as the various subgroups
are encountered during the orbit identification process, they are group
flyped to the position of the first subgroup, thereby transforming the
split group into a full group. At the end of the orbit identification for
the group in question, the size $m$ of its full group will have been
determined. If it is the $n^{th}$ group of size $m$ that has been processed
so far, and there were $k\ge1$ positions in the orbit of the group, then in
an arbitrary manner, we assign the labels $\pm m_n^0$, $\pm m_n^1,\ldots,
\pm m_n^{k-1}$ to the $k$ positions. In the array as constructed so far,
this is recorded as follows. If the two arcs that separate the two
min-tangles which identify the position to which we wish to assign the label
$m_n^i$ are denoted by $e$ and $f$, then we place one copy of the label
$\pm m_n^i$ between the two groups that are the endpoints of $e$ and one
copy of the same label between the two groups that represent the endpoints
of $f$. From this point on, in the identification of the orbits of other
groups still to be analyzed, these two labels of the group's position are
treated as if the group were actually residing at that position (this will
establish the correct order of flyping position on each arc as established
in Theorem \ref{non interference revisited}). In actual fact, one or both
of the ``endpoints'' of $e$ and $f$ may be labels that were deposited
during the identification of the orbit of a group that came earlier in
the list. As we remarked above, once a label is put in place, then for
subsequent group orbit analyses, the label is treated as if it were a
group (whose name is the full label, so it will not be confused with any
other position label for the same group).

The actual procedure to be followed in the identification of the min-tangles
of the orbit depends on whether the group is a loner, a positive
group or a negative group. If a group is a loner, it is first treated as if
it were a negative group and checked to see if it has any flype moves along
this alignment. If so, then it does have an orbit along this alignment
(and therefore not along its positive alignment), and processing continues
to identify the min-tangles of its orbit. On the other hand, if the loner
has no flype moves along its negative alignment, it is then treated as a
positive group. If it has any flype moves along this alignment, then
processing continues to identify the min-tangles of its orbit. If it has
no flype moves along either the negative or the positive alignment, then
it has a trivial orbit.

The decision to check a loner along its negative alignment first (of
course, from a theoretical point of view, the order of checking is not
relevant) was empirically based. Our experience has shown the vast majority
of loners that do have a nontrivial orbit do so along the negative alignment.

It therefore suffices to describe the methods for identifying positive
and negative group orbits.

\begin{figure}[ht]
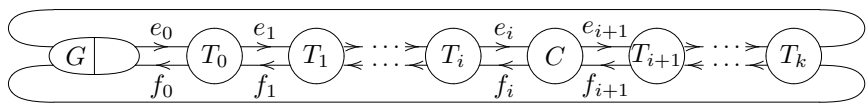

  \centering
  \begin{tabular}{c}
   $\vcenter{\xy /r14pt/:,
      (2,.25)="a",
      {\hloop\hloop-}="a",
      "a";"a"-(0,1)**\dir{-},
      "a"+(-.5,-.5)*{\hbox{$G$}},
      "a"+(3.25,-.5)="t0",{\ellipse(.75){}},
      "t0"*{\hbox{$T_0$}},
      "t0"+(2.75,0)="t1",
      {\ellipse(.75){}},
      "t1"*{\hbox{$T_1$}},
    "a"+(1.1,-.25)="a0","t0"-(.7,-.25)="b0",
    "a0";"b0"**\dir{-}?(.6)*\dir{>}?(.5)*!<0pt,-6pt>{e_0},
    "a"+(1.1,-.75)="a1","t0"-(.7,.25)="b1",
    "b1";"a1"**\dir{-}?(.6)*\dir{>}?(.5)*!<0pt,8pt>{f_0},
    "t0"+(.7,.25)="a2","t1"+(-.7,.25)="b2",
    "a2";"b2"**\dir{-}?(.6)*\dir{>}?(.5)*!<0pt,-6pt>{e_1},
    "t0"+(.7,-.25)="a3","t1"+(-.7,-.25)="b3",
    "b3";"a3"**\dir{-}?(.6)*\dir{>}?(.5)*!<0pt,8pt>{f_1},
    (11.7,-.25)="ti",{\ellipse(.75){}},
    "ti"*{\hbox{$T_i$}},
    "ti"+(-1.3,.25)="c1","ti"+(-.7,.25)="c1a",
    "c1";"c1a"**\dir{-}?(.6)*\dir{>},
    "ti"+(-1.3,-.25)="c2","ti"+(-.7,-.25)="c2a",
    "c2a";"c2"**\dir{-}?(.8)*\dir{>},
    "t1"+(.75,.25)="a5","t1"+(1.25,.25)="b5",
    "a5";"b5"**\dir{-}?(.8)*\dir{>},
     "t1"+(.75,-.25)="a6","t1"+(1.25,-.25)="b6",
     "b6";"a6"**\dir{-}?(.8)*\dir{>},
     "t1"+(1.875,.25)*{\hbox{\small$\ldots$}},
     "t1"+(1.875,-.25)*{\hbox{\small$\ldots$}},
      "ti"+(2.75,0)="core",
      {\ellipse(.75){}},
      "core"*{\hbox{$C$}},
    "ti"+(.7,.25)="c2","core"+(-.7,.25)="d2",
    "c2";"d2"**\dir{-}?(.6)*\dir{>}?(.5)*!<0pt,-6pt>{e_i},
    "ti"+(.7,-.25)="c3","core"+(-.7,-.25)="d3",
    "d3";"c3"**\dir{-}?(.6)*\dir{>}?(.5)*!<0pt,8pt>{f_i},      
      "core"+(2.75,0)="ti+1",
      {\ellipse(.75){}},
      "ti+1"*{\hbox{$T_{i+1}$}},
    "core"+(.7,.25)="ac2","ti+1"+(-.7,.25)="ad2",
    "ac2";"ad2"**\dir{-}?(.6)*\dir{>}?(.5)*!<0pt,-6pt>{e_{i+1}},
    "core"+(.7,-.25)="ac3","ti+1"+(-.7,-.25)="ad3",
    "ad3";"ac3"**\dir{-}?(.6)*\dir{>}?(.5)*!<0pt,8pt>{f_{i+1}},            
    "ti+1"+(3.7,0)="tk",{\ellipse(.75){}},
    "tk"="tx",
    "tk"*{\hbox{$T_{k}$}},
    "tk"+(-1.3,.25)="e1","tk"+(-.7,.25)="e1a",
    "e1";"e1a"**\dir{-}?(.6)*\dir{>},
    "tk"+(-1.3,-.25)="e2","tk"+(-.7,-.25)="e2a",
    "e2a";"e2"**\dir{-}?(.8)*\dir{>},
    "ti+1"+(.75,.25)="aa5","ti+1"+(1.25,.25)="ab5",
    "aa5";"ab5"**\dir{-}?(.8)*\dir{>},
     "ti+1"+(.75,-.25)="aa6","ti+1"+(1.25,-.25)="ab6",
     "ab6";"aa6"**\dir{-}?(.8)*\dir{>},
     "ti+1"+(1.875,.25)*{\hbox{\small$\ldots$}},
     "ti+1"+(1.875,-.25)*{\hbox{\small$\ldots$}},
     "a"+(-1.1,-1.75)="a4","tx"+(.7,-1.25)="b4",
     "a4";"b4"**\dir{-},
     "a"+(-1.1,.75);"tx"+(.7,1.25)**\dir{-},
    "a"+(-1.1,.75),{\hloop-},"a"+(-1.1,-.75),{\hloop-},
     "tx"+(.7,1.25),{\hloop},
     "tx"+(.7,-.25),{\hloop},
    \endxy}$
  \end{tabular}
  \caption{The orbit of a negative group}
\label{orbit of negative group}
\end{figure}

\noindent{\bf Case 1: the orbit of a negative group $G$}. In  Figure
\ref{orbit of negative group}, we have illustrated the general situation.
The core of $G$'s orbit is denoted by $C$ in the figure. The arcs between
two consecutive min-tangles have been labelled so that $e_j$ is the arc
going from left to right and $f_j$ is the arc going the opposite direction.
For the sake of clarity, we have drawn the diagram as if the arc
$e_j$ was above the arc $f_j$, but of course it could be just the opposite.
It also could be that the orbit has only one min-tangle, namely $C$. Of
course, this will come out in the analysis of the orbit.

 If the group code is followed from the first arc of $G$ to the second
 arc of $G$, then we will traverse both arcs of every group belonging
 to any min-tangle of the orbit of $G$ that lies between $G$ and the core of
 the orbit of $G$; namely the tangles $T_0$, $T_1,\ldots,T_i$ as shown in
 Figure \ref{orbit of negative group}. Any group for which exactly one of
 its group arcs is encountered during the traversal from the first arc
 of $G$ through the core and back to the second arc of $G$ must belong to
 the core itself. If we were to continue the traversal from the second
 arc of $G$ through the core and back to the first arc of $G$ (thereby
 completing the knot traversal), all the while making note of the groups
 for which exactly one of the group arcs was encountered in this part of
 the traversal, exactly the same collection of groups would be formed.

 For simplicity, we shall cycle the current array to bring one of the arcs
 of $G$ to the front, and we shall refer to this arc as the first arc of $G$.
 Further, let us refer to the portion of the array
 that lies between the first and the second arcs of $G$ as the first section,
 and the portion which lies between the second arc of $G$ and the end of the
 array (the group arc which precedes the first arc of $G$) as the second
 section.

 With this terminology, the first step in the identification of the orbit of
 $G$ is to determine which groups of the configuration have one arc
 in each section. Having done this, we split the analysis into two parts.

\noindent {\bf Part I}. This part describes the processing of the second
 section. It begins with the identification of the two
 arcs which separate the core from the first
 min-tangle which follows it ($e_{i+1}$ and $f_{i+1}$ in Figure
 \ref{orbit of negative group}). Observe that this is not the min-tangle that
 is encountered after the core has been entered for the first time. The
 knot traversal that the group code provides starts at the first arc of
 $G$  and proceeds through possibly several min-tangles ($T_0$, $T_1,\ldots,
 T_i$) before entering the core $C$. After exiting $C$ for the first time,
 we return to the second arc of $G$ through these same min-tangles (in the
 reverse order), then through possibly more min-tangles ($T_{k}$, $T_{k-1},
 \ldots,T_{i+1}$) until the core is re-entered by means of arc $f_{i+1}$,
 traversed and exited by means of arc $e_{i+1}$, then back through
 $T_{i+1}$, $T_{i+2},\ldots,T_k$ and finally to the first arc of $G$ to
 complete the traversal.

 Once $e_{i+1}$ and $f_{i+1}$ have been identified, the next step will be
 to identify the entry and exit pair which connect $T_{i+1}$ to $T_{i+2}$,
 then the pair which connect $T_{i+2}$ to $T_{i+3}$ and so on, until we
 have worked our way completely through the second section. The last
 pair of arcs that this process identifies is the pair that connects
 the last min-tangle of this portion of the orbit to the group $G$.

 So we must begin by finding the two arcs $e_{i+1}$, $f_{i+1}$ which separate
 the core from the min-tangle that follows the core in the second section.
 We have already identified the groups for which one arc lies in the
 first section and one lies in the second section--these groups we know
 to be in the core, and we shall refer to them as the {\it starter core
 group arcs}. Imagine the strand of the knot which starts at $f_{i+1}$,
 runs through the core, and finishes at $e_{i+1}$. This strand will at
 intervals intertwine with the strand which starts at $e_{i}$, runs
 through the core and stops at $f_{i}$, thereby forming the starter
 core groups. As well, it may intertwine with itself to form additional
 groups. Thus we see that we are searching for the shortest strand which
 contains the starter group arcs and which is closed in the sense that
 any non-starter group which has an arc on this strand has both arcs on
 the strand. To begin with, we extract from the current array the shortest
 subsequence $S_1$ that contains all of the starter group arcs. There
 are two possibilities: either $S_1$ is closed, or else
 there exists at least one group which has exactly one of its arcs in $S_1$.
 In the latter case, we extend $S_1$ to $S_2$, the shortest subsequence
 of the current array that contains the starter arcs and both arcs
 of each group in $S_1$ other than the starter groups. If $S_2$ is not
 closed, repeat this process. Eventually we arrive at a closed subsequence
 $S$ of the current array. Then $f_{i+1}$ is the arc which leads into the
 first group arc of $S$, and $e_{i+1}$ is the arc which leads out
 of the last group of $S$.

 Next, we must find the pair of arcs $e_{i+2}$ and $f_{i+2}$ that separate
 $T_{i+1}$ from $T_{i+2}$ (if indeed there is another min-tangle in this
 segment of the orbit). Form a new set $S_1$ by adjoining to the set $S$
 either the group arc that is adjacent to $e_{i+1}$ but not in the core
 (that is, not in $S$) or the group arc that is adjacent to $f_{i+1}$ but
 not in the core. Again, the selected group arc might be a position marker
 rather than one of the group arcs from the original group code. Close
 $S_1$ as before to obtain a subsequence of the current array. The arcs
 at each end are the sought-after pair $e_{i+2}, f_{i+2}$, unless the
 tangle that lies between the pair $e_{i+1}$ and $f_{i+1}$ and the pair
 just determined is actually another subgroup of $G$ (it is readily verified
 that such a tangle is actually a subgroup of the negative group if and only
 if the tangle itself is just a negative group). In the event that the tangle
 is a subgroup of $G$, we delete this group from the current array and
 from the control list, and the number of crossings it contributes to the
 full group is recorded.

 Continue this process until the subsequence that is formed is the entire
 second section, at which point this part is complete.

\noindent {\bf Part II}.
 This part describes the processing of the first section, whereby 
 the min-tangles which preceed the core are identified. This
 is accomplished by applying the procedure of Part I to the group code
 (or more precisely, the group code as modified by the construction
 procedure so far) that is obtained by writing the group code in the
 reverse order. That is, we construct a temporary working copy of the
 current state of the master group code by starting at the initial group code,
 and write the array in reverse order, so that what was called the first
 section is now the second section. Now apply Part I to this array, but
 place all labels in the appropriate position in the current array, not the
 temporary working copy.

\begin{figure}[ht]
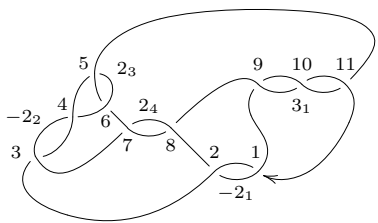

  \centering
  \begin{tabular}{c}
   $\vcenter{   \xy /r33pt/:,
     (3.875,1.625)="a11",
     "a11";(2.875,.625)="a1"*=<4pt,4pt>{\hbox to 2pt{\hfil}}**\crv{(4,1.5) & (4,1.25) & (3.5,.5) & (3,.5)}?(.95)*\dir{}="c"-<-10pt,3pt>;"c"**\dir{}*\dir{>},
    "a1"*=<4pt,4pt>{\hbox to 2pt{\hfil}};(2.375,.625)="a2"**\crv{(2.75,.75) & (2.5,.75)},
    "a2";(.32,.8)="a3"*=<4pt,4pt>{\hbox to 2pt{\hfil}}**\crv{(1.75,0) & (.5,0) & (0,.5)},
    "a3"*=<6pt,6pt>{\hbox to 2pt{\hfil}};(.75,1.25)="a4"**\crv{(.5,.75) & (.75,1)},
    "a4";(1,1.75)="a5"*=<4pt,4pt>{\hbox to 2pt{\hfil}}**\crv{(.75,1.625)},
    "a5"*=<4pt,4pt>{\hbox to 2pt{\hfil}};(1.125,1.375)="a6"**\crv{(1.25,1.75) & (1.25,1.5)},
    "a6";"a4"*=<4pt,4pt>{\hbox to 2pt{\hfil}}**\crv{(1,1.25)},
    "a4"*=<4pt,4pt>{\hbox to 2pt{\hfil}};"a3"**\crv{(.5,1.25) & (.25,1)},
    "a3";(1.375,1.125)="a7"*=<4pt,4pt>{\hbox to 2pt{\hfil}}**\crv{(.5,.5) & (.75,.5)},
    "a7"*=<4pt,4pt>{\hbox to 2pt{\hfil}};(1.875,1.125)="a8"**\crv{(1.5,1.25) & (1.75,1.25)},
    "a8";"a2"*=<4pt,4pt>{\hbox to 2pt{\hfil}}**\dir{-},
    "a2"*=<4pt,4pt>{\hbox to 2pt{\hfil}};"a1"**\crv{(2.5,.5) & (2.75,.5)},
    "a1";(2.875,1.625)="a9"*=<4pt,4pt>{\hbox to 2pt{\hfil}}**\crv{(3,.75) & (3,1) & (2.75,1.25) & (2.75,1.5)},
    "a9"*=<4pt,4pt>{\hbox to 2pt{\hfil}};(3.375,1.625)="a10"**\crv{(3,1.75) & (3.25,1.75)},
    "a10";(3.875,1.625)="a11"*=<4pt,4pt>{\hbox to 2pt{\hfil}}**\crv{(3.5,1.5) & (3.75,1.5)},
    "a11"*=<4pt,4pt>{\hbox to 2pt{\hfil}};"a5"**\crv{(4.25,2) & (4.25,2.25) & (3.75,2.5) & (1.5,2.5) & (1,2)},
    "a5";"a6"*=<4pt,4pt>{\hbox to 2pt{\hfil}}**\crv{(1,1.5)},
    "a6"*=<4pt,4pt>{\hbox to 2pt{\hfil}};"a7"**\dir{-},
    "a7";"a8"*=<4pt,4pt>{\hbox to 2pt{\hfil}}**\crv{(1.5,1) & (1.75,1)},
    "a8"*=<4pt,4pt>{\hbox to 2pt{\hfil}};"a9"**\crv{(2.5,1.75) & (2.75,1.75)},
    "a9";"a10"*=<4pt,4pt>{\hbox to 2pt{\hfil}}**\crv{(3,1.5) & (3.25,1.5)},
    "a10"*=<4pt,4pt>{\hbox to 2pt{\hfil}};"a11"**\crv{(3.5,1.75) & (3.75,1.75)},
            (2.625,.375)*{\hbox{$\ssize -2_1$}},
            (.19,1.25)*{\hbox{$\ssize -2_2$}},
            (1.375,1.8)*{\hbox{$\ssize 2_3$}},
            (1.625,1.375)*{\hbox{$\ssize 2_4$}},
            (3.375,1.375)*{\hbox{$\ssize 3_1$}},
            (2.85,.85)*{\hbox{$\ssize 1$}},
            (2.375,.85)*{\hbox{$\ssize 2$}},
            (.1,.85)*{\hbox{$\ssize 3$}},
            (.63,1.4)*{\hbox{$\ssize 4$}},
            (.875,1.87)*{\hbox{$\ssize 5$}},
            (1.125,1.2)*{\hbox{$\ssize 6$}},
            (1.375,.925)*{\hbox{$\ssize 7$}},
            (1.875,.925)*{\hbox{$\ssize 8$}},
            (2.875,1.86)*{\hbox{$\ssize 9$}},
            (3.375,1.86)*{\hbox{$\ssize 10$}},
            (3.875,1.86)*{\hbox{$\ssize 11$}},
      \endxy}$
  \end{tabular}
  \caption{A configuration of a prime 11 crossing knot}
\label{example for master group code construction}
\end{figure}

 We illustrate this process with an example. Consider the knot configuration
 shown in Figure \ref{example for master group code construction}. This is a
 prime alternating 11 crossing knot in full group. We have shown a
 Gauss labelling on the diagram, and if we start at the crossing labelled
 1 and proceed in the direction of the crossing labelled 2 on the strand
 taking us to crossing labelled 3, we obtain the following group code for
 the configuration.
 $$
  -2_1,-2_2,2_3,-2_2,2_4,-2_1,3_1,2_3,2_4,3_1
 $$ 
 
 The control list of groups is $-2_1,-2_2,2_3,2_4,3_1$, and so we would
 begin by finding the orbit of the negative group $G=-2_1$. No cycling
 of the array is necessary, since one of the group arcs of $-2_1$ is
 the first entry in the array. The first and second sections, respectively,
 are
 $$
  -2_2,2_3,-2_2,2_4  \quad\hbox{and}\quad 3_1,2_3,2_4,3_1
 $$
 and so the starter group arcs are $2_3$ and $2_4$. Since the
 sequence $S_1=2_3,2_4$ is closed, $S=S_1$ and we have found the
 entry and exit arcs for this side of the core of the orbit of $G$
 to be the arc from 11 to 5 and the arc from 8 to 9, respectively. 
 In this particular example, the two adjacent arcs to the subsequence $S$
 belong to the same group, $3_1$. We form $S_1=3_1,2_3,2_4$. This is not
 closed, since the second arc of group $3_1$ is missing. Upon closing
 it up, we obtain $S=3_1,2_3,2_4,3_1$. The tangle which has just been
 determined has incident arcs $(5,11)$ and $(8,9)$ as determined earlier,
 and the arcs $(1,9)$ and $(1,11)$ which have just been determined. This
 tangle consists of just the group $3_1$, which must therefore be checked
 to see if it is a subgroup of the full group determined by $-2_1$. Since
 it is a positive group, it is not a subgroup of the negative group to
 which $-2_1$ belongs, so it is a min-tangle in the orbit of $-2_1$.
 Finally, since $S$ is equal to the second section, part I has been
 completed.

 For part II, we reverse the array to obtain
 $$
   -2_1,3_1,2_4,2_3,3_1,-2_1,2_4,-2_2,2_3,-2_2
 $$
 with first and second sections
 $$
    3_1,2_4,2_3,3_1 \quad\hbox{and}\quad 2_4,-2_2,2_3,-2_2,
 $$
 respectively. We know that the starter arcs are $2_3$ and $2_4$,
 so $S_1=2_4,-2_2,2_3$. This is not closed, since the other arc for
 group $-2_2$ is missing. Upon closing $S_1$, we obtain $S=2_4,-2_2,
 2_3,-2_2$, which is the second section. Thus the entry and exit arcs
 for this side of the core of $G=-2_1$ are $(2,3)$ and $(8,2)$,
 respectively. Since we have exhausted the second section (the first
 section of the actual array, since we are doing part II), there are
 no additional min-tangles. Thus the orbit of $G=-2_1$ is of the
 form

 \begin{figure}[ht]
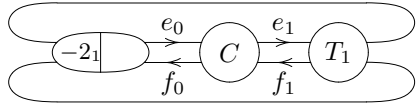

  \centering
  \begin{tabular}{c}
   $\vcenter{\xy /r15pt/:,
      (2,.25)="a",
      {\hloop\hloop-}="a",
      "a";"a"-(0,1)**\dir{-},
      "a"+(-.5,-.5)*{\hbox{\small $-2_1$}},
      "a"+(3.25,-.5)="t0",{\ellipse(.75){}},
      "t0"*{\hbox{$C$}},
      "t0"+(2.75,0)="t1",
      {\ellipse(.75){}},
      "t1"*{\hbox{$T_1$}},
    "a"+(1.1,-.25)="a0","t0"-(.7,-.25)="b0",
    "a0";"b0"**\dir{-}?(.6)*\dir{>}?(.5)*!<0pt,-6pt>{e_0},
    "a"+(1.1,-.75)="a1","t0"-(.7,.25)="b1",
    "b1";"a1"**\dir{-}?(.6)*\dir{>}?(.5)*!<0pt,8pt>{f_0},
    "t0"+(.7,.25)="a2","t1"+(-.7,.25)="b2",
    "a2";"b2"**\dir{-}?(.6)*\dir{>}?(.5)*!<0pt,-6pt>{e_1},
    "t0"+(.7,-.25)="a3","t1"+(-.7,-.25)="b3",
    "b3";"a3"**\dir{-}?(.6)*\dir{>}?(.5)*!<0pt,8pt>{f_1},
    "t1"="tx",
     "a"+(-1.1,-1.75)="a4","tx"+(.7,-1.25)="b4",
     "a4";"b4"**\dir{-},
     "a"+(-1.1,.75);"tx"+(.7,1.25)**\dir{-},
    "a"+(-1.1,.75),{\hloop-},"a"+(-1.1,-.75),{\hloop-},
     "tx"+(.7,1.25),{\hloop},
     "tx"+(.7,-.25),{\hloop},
    \endxy}$
  \end{tabular}
  \caption{The orbit of $-2_1$}
\label{orbit of negative 21}
\end{figure}

\noindent and the updated array, showing the additional position of group
$-2_1$ is
$$
  -2_1^0,-2_2,2_3,-2_2,2_4,-2_1^0,3_1,-2_1^1,2_3,2_4,-2_1^1,3_1.
$$
\noindent Now that the orbit of $-2_1$ has been identified, we can remove
it from the control list. The revised control list is $-2_2,2_3,2_4,3_1$, and
the next group in the list is also a negative group. We cycle the current
array so as to start with one of the group arcs for $G=-2_2$, say
$$
  -2_2,2_3,-2_2,2_4,-2_1^0,3_1,-2_1^1,2_3,2_4,-2_1^1,3_1,-2_1^0.
$$
The first and second sections are $2_3$ and $2_4,-2_1^0,3_1,-2_1^1,2_3,2_4,-2_1^1,3_1,-2_1^0$,
respectively, and there is only one starter group arc; namely $2_3$. We
have $S_1=2_3=S$, and the entry and exit arcs from this side of the core
are $(11,5)$ and $(6,7)$, respectively. The adjacent group arcs are
$-2_1^1$ and $2_4$, so we choose one, say $-2_1^1$, and form $S_1=-2_1^1,
2_3$. Since the other arc of group $-2_1^1$ is missing, $S_1$ is not closed.
We find the other arc of group $-2_1^1$ and form $S_2=-2_1^1,2_3,2_4,-2_1^1$.
Since the other arc of group $2_4$ is missing, we find it and form
$S_3=2_4,-2_1^0,3_1,-2_1^1,2_3,2_4,-2_1^1$. Since $S_3$ is missing the other
arc from each of groups $-2_1^0$ and $3_1$, we form $S_4=
2_4,-2_1^0,3_1,-2_1^1,2_3,2_4,-2_1^1,3_1,-2_1^0$, which is closed. Thus
$S=S_4$ is the second section, so this part is done. For part II, we
reverse the sequence to obtain
$$
  -2_2,-2_1^0,3_1,-2_1^1,2_4,2_3,-2_1^1,3_1,-2_1^0,2_4,-2_2,2_3,
$$
whose second section is simply $2_3$. Thus there is nothing to do in
part II. The core of the orbit of group $-2_2$ is just the group $2_3$.
The updated array, showing the additional position of group $-2_2$ is
$$
  -2_2^0,2_3,-2_2^0,2_4,-2_1^0,3_1,-2_1^1,-2_2^1,2_3,-2_2^1,2_4,-2_1^1,
  3_1,-2_1^0.
$$
We delete the group $-2_2$ from the control list, thereby obtaining the
revised control list $2_3,2_4,3_1$. Since the remaining groups are positive,
we must wait until after the method for identifying the orbit of a positive
group has been discussed before completing the construction of the master
group code for this example.

\noindent{\bf Case 2: the orbit of a positive group $G$}. The situation
for a positive group is considerably simpler than that of a negative
group. As before, we cycle the current array so as to begin with one
of the arcs of $G$, which we shall refer to as the first arc of $G$, and
we refer to the portion of the array
that lies between the first and the second arcs of $G$ as the first section,
and the portion which lies between the second arc of $G$ and the end of the
array as the second section. The arcs labelled $e_i$ join group arcs
that are in the first section, while the arcs labelled $f_i$ join group
arcs that are in the second section (again, we have taken some artistic liberties,
and listed $e_i$ as if it was always above $f_i$ in the configuration, but of
course this is not necessarily the case).

 \begin{figure}[ht]
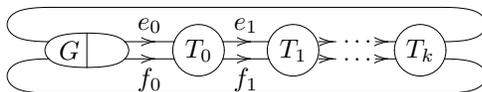

  \centering
  \begin{tabular}{c}
   $\vcenter{\xy /r13pt/:,
      (2,.25)="a",
      {\hloop\hloop-}="a",
      "a";"a"-(0,1)**\dir{-},
      "a"+(-.5,-.5)*{\hbox{$G$}},
      "a"+(3.25,-.5)="t0",{\ellipse(.75){}},
      "t0"*{\hbox{$T_0$}},
      "t0"+(2.75,0)="t1",
      {\ellipse(.75){}},
      "t1"*{\hbox{$T_1$}},
    "a"+(1.1,-.25)="a0","t0"-(.7,-.25)="b0",
    "a0";"b0"**\dir{-}?(.6)*\dir{>}?(.5)*!<0pt,-6pt>{e_0},
    "a"+(1.1,-.75)="a1","t0"-(.7,.25)="b1",
    "a1";"b1"**\dir{-}?(.6)*\dir{>}?(.5)*!<0pt,8pt>{f_0},
    "t0"+(.7,.25)="a2","t1"+(-.7,.25)="b2",
    "a2";"b2"**\dir{-}?(.6)*\dir{>}?(.5)*!<0pt,-6pt>{e_1},
    "t0"+(.7,-.25)="a3","t1"+(-.7,-.25)="b3",
    "a3";"b3"**\dir{-}?(.6)*\dir{>}?(.5)*!<0pt,8pt>{f_1},
    "t1"+(.75,.25)="a5","t1"+(1.25,.25)="b5",
    "a5";"b5"**\dir{-}?(.8)*\dir{>},
     "t1"+(.75,-.25)="a6","t1"+(1.25,-.25)="b6",
     "a6";"b6"**\dir{-}?(.8)*\dir{>},
     "t1"+(1.875,.25)*{\hbox{\small$\ldots$}},
     "t1"+(1.875,-.25)*{\hbox{\small$\ldots$}},
    "t1"+(3.7,0)="tk",{\ellipse(.75){}},
    "tk"+(-1.3,.25)="c1","tk"+(-.7,.25)="c1a",
    "c1";"c1a"**\dir{-}?(.6)*\dir{>},
    "tk"+(-1.3,-.25)="c2","tk"+(-.7,-.25)="c2a",
    "c2";"c2a"**\dir{-}?(.6)*\dir{>},    
    "tk"="tx",
    "tk"*{\hbox{$T_{k}$}},
     "a"+(-1.1,-1.75)="a4","tx"+(.7,-1.25)="b4",
     "a4";"b4"**\dir{-},
     "a"+(-1.1,.75);"tx"+(.7,1.25)**\dir{-},
    "a"+(-1.1,.75),{\hloop-},"a"+(-1.1,-.75),{\hloop-},
     "tx"+(.7,1.25),{\hloop},
     "tx"+(.7,-.25),{\hloop},
    \endxy}$
  \end{tabular}
  \caption{The orbit of a positive group}
\label{orbit of positive group}
\end{figure}

Suppose that we have identified the arcs $e_i$ and $f_i$ that identify
the $i^{th}$ position of the group $G$ (that is, the two arcs that
separate min-tangle $T_i$ from min-tangle $T_{i+1}$, if $i>0$, or the
pair of arcs $e_0$ and $f_0$ that connect $G$ to $T_0$). Let $S_i^1$ denote
the sequence of all group arcs that belong to the min-tangles between $G$ and
the position determined by $e_i$ and $f_i$ in the order that they are
encountered on the strand which leaves $G$ along $e_0$, and let $S_i^2$ denote
the sequence of all group arcs that belong to the same min-tangles but are
encountered on the strand which leaves $G$ along $f_0$. Thus $S_i^1$ is an
initial segment of the first section, and $S_i^2$ is an initial segment
of the second section, and these two segments constitute a closed pair,
by which we mean that any group that has at least one of its arcs appear
in either $S_i^1$ or
$S_i^2$ also has the other arc appearing in either $S_i^1$ or $S_i^2$.
Note that $S_0^1$ and $S_0^2$ are empty. If either of $S_i^1$ or $S_i^2$
is equal to the section to which it belongs, the orbit has been completely
identified. Let $S_{i+1}^1$ and $S_{i+1}^2$ be the shortest
initial segments of the first and second sections, respectively, for which
$S_i^1$ is a proper initial segment of $S_{i+1}^1$ and $S_{i}^2$ is a
proper initial segment of $S_{i+1}^2$ which are a closed pair in the sense described
above. The groups whose arcs appear in $S_{i+1}^1$ or $S_{i+1}^2$ but not
in $S_i^1$ or $S_i^2$ form a tangle which is either a min-tangle in the
orbit of $G$ or a subgroup of the full group determined by $G$. It is
easy to verify that the tangle is a subgroup of $G$ if and only if the
tangle consists of a single positive group, in which case its two group
arcs are deleted from the current array and from the first and second
sections, it is removed from the control list, the number of crossings
it contributes to the size of the full group to which $G$ belongs is
recorded, and we reset $S_{i+1}^1=S_i^1$ and $S_{i+1}^2=S_i^2$, and
repeat the process (in effect, we have group flyped the subgroup to
amalgamate it with the original group $G$). On the other hand, if the
tangle is a min-tangle in the orbit of $G$, then we record the pair of
arcs $e_i$ and $f_i$ as a position pair for the orbit of $G$, and set
$e_{i+1}$ to be the arc joining the last group arc of $S_{i+1}^1$ to
the next group arc to be encountered on that strand (either the next
group arc in the first section, or, if the first section has been
exhausted, to the arc of $G$ that follows the first section), and
set $f_{i+1}$ to be the arc joining the last group arc of $S_{i+1}^2$ to
the next group arc to be encountered on that strand (either the next
group arc in the second section, or, if the second section has been
exhausted, to the arc of $G$ that follows the second section), and
repeat the process, replacing $i$ by $i+1$.

When the identification of the orbit is complete, the size of the
full group, say $m$, is known, and the full group is assigned the index of
the next group of size $m$ to be labelled, say $n$. Suppose that there are
$k$ positions in the orbit of $G$. The current array is then modified as
follows: the two arcs of $G$ are
replaced by the label $m_n^0$, and if $k>1$, then for each arc in each pair
$e_i$, $f_i$, $i=1,\ldots, k-1$, the label $m_n^i$ is inserted between the
two group arcs that are joined by the arc in question. This completes the
processing of the orbit of $G$.

We illustrate this procedure by completing the construction of the master
group code for the knot shown in Figure \ref{example for master group code construction}.
We had left off with the control list $2_3,2_4,3_1$, so the positive group
$2_3$ is the next to be processed. We cycle the current array so as to
begin with an arc of group $2_3$. The result is:
$$
  2_3,-2_2^0,2_4,-2_1^0,3_1,-2_1^1,-2_2^1,2_3,-2_2^1,2_4,-2_1^1,3_1,-2_1^0,
  -2_2^0
$$
and so $-2_2^0,2_4,-2_1^0,3_1,-2_1^1,-2_2^1$ and $-2_2^1,2_4,-2_1^1,3_1,
-2_1^0,-2_2^0$ are the first and second sections, respectively. We start with
$S_0^1$ and $S_0^2$ empty, and look for the shortest initial subsequences
of the first and second sections that contain $S_0^1$ and $S_0^2$, respectively,
and which form a closed pair. We find that $S_1^1$ is the first section and
$S_1^2$ is the second section, so the orbit of $2_3$ has a single min-tangle.
We replace both occurrences of $2_3$ by $2_3^0$, and remove group $2_3$ from
the control list to complete the processing of the orbit of group $2_3$. The
control list is now $2_4,3_1$, so we cycle the current array so as to begin
with an arc of group $2_4$, say
$$
  2_4,-2_1^0,3_1,-2_1^1,-2_2^1,2_3^0,-2_2^1,2_4,-2_1^1,3_1,-2_1^0,
  -2_2^0,2_3^0,-2_2^0,
 $$
 and we find that $-2_1^0,3_1,-2_1^1,-2_2^1,2_3^0,-2_2^1$ and $-2_1^1,3_1,
 -2_1^0,-2_2^0,2_3^0,-2_2^0$ are the first and second sections. We start with
 $S_0^1$ and $S_0^2$ empty and look for the shortest initial subsequences
 $S_1^1$ and $S_1^2$ of the first and second sections which properly contain
 $S_0^1$ and $S_0^2$, respectively, and which form a closed pair. We find that
 $S_1^1=-2_1^0,3_1,-2_1^1$ and $S_1^2=-2_1^1,3_1,-2_1^0$, so $e_1$ is the arc
 that joins $-2_1^1$ to $-2_2^1$ in the first section, and $f_1$ is the arc
 that joins $2_1^0$ to $-2_2^0$ in the second section. At the next step, we
 find that $S_2^1$ is the first section, so $S_2^2$ must be the second
 section, and the orbit has been completely identified. There are two
 min-tangles in the orbit of $2_4$. We replace both occurrences of $2_4$ by
 $2_4^0$, and place the label $2_4^1$ between the first occurrences of $-2_1^1$
 and $-2_2^1$, and between the second occurrences of $-2_1^0$ and $-2_2^0$. The
 result is
 $$
  2_4^0,-2_1^0,3_1,-2_1^1,2_4^1,-2_2^1,2_3^0,-2_2^1,2_4^0,
                       -2_1^1,3_1,-2_1^0,2_4^1,-2_2^0,2_3^0,-2_2^0.
 $$
 We may now remove group $2_4$ from
 the control list to complete the processing of the orbit of group $2_4$.
 The control list is now $3_1$, so we cycle the current array so as to begin
 with an arc of group $3_1$, say
 $$
  3_1,-2_1^1,2_4^1,-2_2^1,2_3^0,-2_2^1,2_4^0,-2_1^1,3_1,-2_1^0,
                             2_4^1,-2_2^0,2_3^0,-2_2^0,2_4^0,-2_1^0
 $$
 with $-2_1^1,2_4^1,-2_2^1,2_3^0,-2_2^1,2_4^0,-2_1^1$ and
 $-2_1^0,2_4^1,-2_2^0,2_3^0,-2_2^0,2_4^0,-2_1^0$ the first and second
 sections, respectively. We start with $S_0^1$ and $S_0^2$ empty and look
 for the shortest initial subsequences $S_1^1$ and $S_1^2$ of the first
 and second sections which properly contain $S_0^1$ and $S_0^2$,
 respectively, and which form a closed pair. Since the first arc of
 the second section is $-2_1^0$, and the second occurrence of
 that group is at the end of the second section, we see that $S_1^1$ and
 $S_1^2$ must be the first and second sections, respectively, and the
 processing of the orbit of group $3_1$ is complete--its orbit has a single
 min-tangle. We replace both occurrences of $3_1$ by $3_1^0$, and remove
 group $3_1$ from the control list to complete the processing of the orbit
 of group $3_1$. Since the control list is now empty, the construction of
 the master group code is complete. The master group code that we have
 constructed from the group code
 $$
  -2_1,-2_2,2_3,-2_2,2_4,-2_1,3_1,2_3,2_4,3_1
 $$ 
 is
 $$
  3_1^0,-2_1^1,2_4^1,-2_2^1,2_3^0,-2_2^1,2_4^0,-2_1^1,
             3_1^0,-2_1^0,2_4^1,-2_2^0,2_3^0,-2_2^0,2_4^0,-2_1^0.
 $$

 Once we have a master group code for a full group knot configuration,
 then by selecting one position for each group and then writing out
 the selected group arcs in the order that they appear in the master
 group code, we will produce a group code for a full group configuration
 of the knot, and every full group configuration of the knot can be
 obtained in this way.

\begin{figure}[ht]
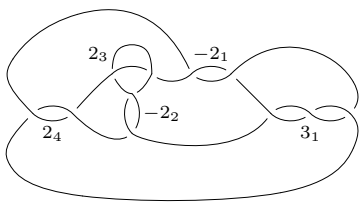

  \centering
  \begin{tabular}{c}
   $\vcenter{\xy /r30pt/:, (4.375,1.125)="a1", "a1";(.375,1.125)="a2"*=<4pt,4pt>{\hbox to 2pt{\hfil}}**\crv{(4.5,1) & (4.5,.75) & (4,0) & (.5,0) & (0,.5) & (0,.75)}, "a2"*=<4pt,4pt>{\hbox to 2pt{\hfil}};(.875,1.125)="a3"**\crv{(.5,1.25) & (.75,1.25)}, "a3";(1.625,.875)="a4"*=<4pt,4pt>{\hbox to 2pt{\hfil}}**\crv{(1.25,.75) & (1.5,.75)}, "a4"*=<4pt,4pt>{\hbox to 2pt{\hfil}};(1.625,1.375)="a5"**\crv{(1.75,1) & (1.75,1.25)},
  "a5";(1.375,1.625)="a6"*=<4pt,4pt>{\hbox to 2pt{\hfil}}**\crv{(1.5,1.375)},
  "a6"*=<4pt,4pt>{\hbox to 2pt{\hfil}};(1.875,1.625)="a7"**\crv{(1.375,1.875) & (1.5,2) & (1.75,2) & (1.875,1.875)},
  "a7";"a5"*=<4pt,4pt>{\hbox to 2pt{\hfil}}**\crv{(1.75,1.375)},
  "a5"*=<4pt,4pt>{\hbox to 2pt{\hfil}};"a4"**\crv{(1.5,1.25) & (1.5,1)},
  "a4";(3.375,1.125)="a8"*=<4pt,4pt>{\hbox to 2pt{\hfil}}**\crv{(1.75,.75) & (3,.5)},
  "a8"*=<4pt,4pt>{\hbox to 2pt{\hfil}};(3.875,1.125)="a9"**\crv{(3.5,1.25) & (3.75,1.25)},
  "a9";"a1"*=<4pt,4pt>{\hbox to 2pt{\hfil}}**\crv{(4,1) & (4.25,1)},
  "a1"*=<4pt,4pt>{\hbox to 2pt{\hfil}};(2.875,1.625)="a10"**\crv{(4.5,1.25) & (4.5,1.5) & (4,2) & (3.25,2)},
  "a10";(2.375,1.625)="a11"*=<4pt,4pt>{\hbox to 2pt{\hfil}}**\crv{(2.75,1.5) & (2.5,1.5)},
  "a11"*=<4pt,4pt>{\hbox to 2pt{\hfil}};"a2"**\crv{(2,2.5) & (.75,2.5) & (0,1.75) & (0,1.5)},
  "a2";"a3"*=<4pt,4pt>{\hbox to 2pt{\hfil}}**\crv{(.5,1) & (.75,1)},
  "a3"*=<4pt,4pt>{\hbox to 2pt{\hfil}};"a6"**\crv{(1.2,1.5)},
  "a6";"a7"*=<4pt,4pt>{\hbox to 2pt{\hfil}}**\crv{(1.5,1.75) & (1.75,1.75)},
  "a7"*=<4pt,4pt>{\hbox to 2pt{\hfil}};"a11"**\crv{(2,1.5) & (2.25,1.5)},
  "a11";"a10"*=<4pt,4pt>{\hbox to 2pt{\hfil}}**\crv{(2.5,1.75) & (2.75,1.75)},
  "a10"*=<4pt,4pt>{\hbox to 2pt{\hfil}};"a8"**\dir{-},
  "a8";"a9"*=<4pt,4pt>{\hbox to 2pt{\hfil}}**\crv{(3.5,1) & (3.75,1)},
  "a9"*=<4pt,4pt>{\hbox to 2pt{\hfil}};"a1"**\crv{(4,1.25) & (4.25,1.25)},
  (2.625,1.875)*{\hbox{$\ssize -2_1$}},
  (2,1.125)*{\hbox{$\ssize -2_2$}},
  (1.2,1.875)*{\hbox{$\ssize 2_3$}},
  (.625,.875)*{\hbox{$\ssize 2_4$}},
  (3.875,.875)*{\hbox{$\ssize 3_1$}},
  \endxy}$
  \end{tabular}
  \caption{A new configuration for the prime 11 crossing knot}
 \label{knot for constructed code}
 \end{figure}

 For example, if we select the zero position for groups $3_1$, $-2_2$ and
 $2_3$, and position one for each of $-2_1$ and $2_4$, we obtain the group
 code
 $$
  3_1,-2_1,2_4,2_3,-2_1,3_1,2_4,-2_2,2_3,-2_2.
 $$
 We have shown this knot in Figure \ref{knot for constructed code}.

\subsection{The Master Array for a Knot}

 Each master group code for a knot provides the orbit structure for the
 knot, detailing for each group its size, sign (positive or negative)
 and the min-tangles of its orbit, and any two master group codes for
 a knot have the same structure, in the sense that any one can be
 obtained from any other by an appropriate relabelling and possibly a
 reversal of direction of traversal.

 However, comparing two master group codes to determine whether or
 not they represent the same knot is a computationally demanding task.
 As an alternative, we introduce a procedure to produce an array
 from a master group code, where this array has the property that
 two master group codes represent configurations of the same knot if
 and only if the constructed arrays for the two master group codes are
 identical. This array shall be called the master array for the knot.

 As a first step, we intoduce a linear ordering on the labels in a master
 group code. Recall that each label is an ordered triple $(m,n,i)$, where
 $m$ is a signed integer whose magnitude is the size of the group that the
 label represents, $n$ is a positive integer which is the index of this
 particular group of size $|m|$, and $i$ is a non-negative integer which
 represents a position in the orbit of the group. We shall use the usual
 ordering of integers to lexicographically order these triples.

 Now, let $M$ be a master group code for a knot $K$. Form the set
 of all sequences that can be obtained by selecting an entry in $M$ and
 a direction to traverse $M$ (recall that $M$ is a cyclic code), then
 record the entries of $M$ as they are encountered when the cycle is followed
 in the selected direction from the selected starting point. Working with
 this set corresponds to rotating (and possibly reversing) a master group
 code in order to compare it to some other master group code. The problem
 of the arbitrariness of the labelling still remains. To deal with this,
 we relabel each sequence in our set. For each sequence, we work from the
 beginning, changing each label $\pm m_n^i$ as it is encountered in the
 sequence according to the following scheme.  We leave $m$ unchanged (since
 this is the size of the group), but $n$ and $i$ will be changed to $p$
 and $j$, respectively, where this label represents the $j+1$ orbit
 position of the $p^{th}$ group of size $m$ that has been encountered by the time
 this occurrence of the label $\pm m_n^i$ has been reached. Thus after
 relabelling, all groups of a given size $m$ will now appear in the sequence
 in the order $m_1, m_2,$ and so on, and the positions of the group labelled
 $m_j$ will now be labelled in order of appearance as $m_j^0$, $m_j^1$, and 
 so on. 

 Now, since each sequence in this set consists of labels $\pm m_n^i$, which
 we may consider as ordered triples $(\pm m,n,i)$ of integers, we could
 lexicographically order the set of sequences that we have obtained from
 the master group code (where we consider the triples to be linearly
 ordered by the lexicographical ordering on $\mathbb{Z}^3$ that is obtained
 from the usual ordering on $\mathbb{Z}$). Then (of many possible choices), we
 might choose the unique minimum sequence as our designated master array.
 While this would be reasonable in principle, the set of sequences that we
 must analyze is large enough to present computational headaches.
 Consequently, we introduce a scheme for selecting a special subset of the
 set of all sequences that we obtain from a master group code, and then
 apply the lexicographical ordering to select the
 least element of this special subset. This least sequence is the one that
 we call the master array for the knot. Since different master group codes
 for the same knot will all produce the same set of sequences, it follows
 that this scheme will result in equal master arrays when applied to
 two master group codes for the same knot. Morever, since the master array
 is a sequence which, when converted to a cycle by having the first element
 of the sequence follow the last element of the sequence, does yield a
 master group code for the knot, it follows that different master arrays
 correspond to different knots.

 \begin{definition}
  Let $K$ be a prime alternating knot and let $M$ be a master group code
  for $K$. A {\it ladder} of $M$ is an interval $m_i,m_{i+1},\ldots,m_{i+k}$
  of consecutive entries of $M$ (where the indices are to be calculated
  modulo the length of $M$ to account for the fact that $M$ is a cyclic
  code) for which no two entries in the interval are of the form
  $(\pm m,n,i)$ and $(\pm m,n,j)$ for any values of $i$ and $j$. The
  {\it length} of a ladder is the sum of the sizes of the groups
  that appear in the ladder.
  Finally, the {\it maximum ladder length} of $K$, denoted by $\lmax{K}$,
  or simply $\lmaxs$ when $K$ is understood, is defined to be
  the maximum of the lengths of the ladders of $M$.
 \end{definition}
 
 As is implied in the definition of $\lmaxs$, it is independent of the
 particular master group code from which it is computed.

 Now, of all of the relabelled sequences that we have formed from the
 master group code that we are working with, we select only those
 that begin with an $\lmaxs$ ladder; that is, a ladder of length
 $\lmaxs$. This set of sequences shall be denoted by $\Lmax$.
 
 \begin{definition}\label{the master array}
  The {\it master array} for a given knot $K$ is the least element
  in the lexicographically ordered set $\Lmax$.
 \end{definition}

 We shall finish up this section by illustrating the construction of
 the master array from the master group code
 $$
  3_1^0,-2_1^1,2_4^1,-2_2^1,2_3^0,-2_2^1,2_4^0,-2_1^1,3_1^0,
                       -2_1^0,2_4^1,-2_2^0,2_3^0,-2_2^0,2_4^0,-2_1^0.
 $$
 for the knot shown in Figure \ref{example for master group code
 construction}. The first job is to determine the ladders of greatest length.
 Conceptually, for each entry in the master group code, we determine the
 longest ladder that begins at the selected entry (following a clockwise
 direction--the direction proceeding from left to right in the linear
 representation of the cycle as shown above). From these, we determine
 those of maximum length, and then for each ladder of maximum length, we
 form two sequences, one obtained by starting at the beginning of the ladder
 and proceeding in the clockwise direction around the cycle, the other
 obtained by starting at the end of the ladder and proceeding in the
 counterclock-wise direction around the cycle. The resulting set of
 sequences, after relabelling, is $\Lmax$. We then seek the least element
 of this set.

 Starting with the element $3_1^0$, we find the ladder $3_1^0,-2_1^1,2_4^1,
 -2_2^1,2_3^0$ which can't be extended further, since the next entry is $-2_2^1$
 and the group $-2_2$ has already made an appearance. The length of this ladder
 is $3+2+2+2+2=11$. Since we are searching for ladders of greatest length, we can
 skip every ladder which begins at any point between $3_1^0$ and the first appearance
 of the label $-2_2^1$, since they will necessarily be shorter. Thus the next
 starting point that we check will be the entry immediately after the first appearance
 of $-2_2^1$, namely $2_3^0$. We find the 
 ladder $2_3^0,-2_2^1,2_4^0,-2_1^1,3_1^0$ which can't be extended, since the next
 entry is $-2_1^0$, and $-2_1^1$ already appears in the ladder. This ladder also has
 length 11. We move ahead to the entry that follows $-2_1^1$, and find the ladder
 $3_1^0,-2_1^0,2_4^1,-2_2^0,2_3^0$, which has length 11. We may again skip ahead to
 the start point $2_3^0$, and there we find the ladder $2_3^0,-2_2^0,2_4^0,-2_1^0, 
 3_1^0$, a ladder of length 11 as well. Since the group that stops this ladder is
 $-2_1$, we skip ahead to the first entry that follows $-2_1^0$. Since this is the
 entry $3_1^0$ which is the first entry that we considered, we have determined that
 $\lmaxs=11$ for this knot. Morever, we have found 4 maximum length ladders in the
 master group code. Since we are using lexicographical ordering, a group of 3 will
 lose out to a group of 2, so we don't need to expend the effort to relabel the 
 sequences that begin with a group of 3. Each of the four ladders of length 11
 produces one sequence that starts with a group of 2 and sequence that starts
 with a group of 3. After discarding the four sequences that begin with a
 group of 3, we are left with 4 sequences that begin with a group of 2.
 The one that we obtain from the first $\lmaxs$ ladder that was found above
 is (before relabelling)
 $$
   2_3^0,-2_2^1,2_4^1,-2_1^1,3_1^0,-2_1^0,2_4^0,-2_2^0,2_3^0,
              -2_2^0,2_4^1,-2_1^0,3_1^0,-2_1^1,2_4^0,-2_2^1,
 $$
 which, after relabelling, becomes
 $$
  2_1^0,-2_2^0,2_3^0,-2_4^0,3_1^0,-2_4^1,2_3^1,-2_2^1,2_1^0,
                       -2_2^1,2_3^0,-2_4^1,3_1^0,-2_4^0,2_3^1,-2_2^0.
 $$
 The next $\lmaxs$ ladder provides the sequence (before relabelling)
 $$
  2_3^0,-2_2^1,2_4^0,-2_1^1,3_1^0,-2_1^0,2_4^1,-2_2^0,2_3^0,
              -2_2^0,2_4^0,-2_1^0,3_1^0,-2_1^1,2_4^1,-2_2^1,
 $$
 which, after relabelling, results in the same sequence as we obtained above.
 The remaining two $\lmaxs$ ladders also produce this sequence after
 relabelling, so the master array for the 11 crossing knot shown in
 Figure \ref{example for master group code construction}
 (and of course,
 for the knot shown in Figure \ref{knot for constructed code}) is
 $$
  2_1^0,-2_2^0,2_3^0,-2_4^0,3_1^0,-2_4^1,2_3^1,-2_2^1,2_1^0,
                       -2_2^1,2_3^0,-2_4^1,3_1^0,-2_4^0,2_3^1,-2_2^0.
 $$

%% file: implem1.tex
\section{An Implementation of the Construction, Part I: Full Group
  Configurations Suffice}

 In many of the various attempts to enumerate alternating prime knots over
 the years, a major stumbling block was the large amount of redundant
 work that had to be carried out. A primary feature of our knot encoding
 scheme allows us to work with a given knot in a full group configuration
 only. By not having to work with split group configurations, a tremendous
 amount of unnecessary work is avoided.

 In this section, we prove that it suffices to apply the operators introduced
 above to knots in full group configuration only, and in the case of $D$,
 to a single crossing of each group that is eligible for $D$ to act on.
 Furthermore, we shall examine each operator in order to identify situations
 that we can be assured will result in redundant constructions. This will
 enable us to instruct each operator to bypass certain situations, thereby
 obtaining a more efficient construction process.

 For conceptual purposes, we quickly review the basic construction process
 here. The construction
 begins with the application of $D$ and $ROTS$ to the knot configurations
 of $n$ crossings. That is, for each negative group, each loner and each
 positive 2-group of any knot configuration of $n$ crossings, $D$ is applied
 to any one crossing of the selected group.

 After all possible applications of $D$ have been made, the next step
 is to apply the $ROTS$ operator to each negative group size 2 or 3 of
 any knot configuration of $n$ crossings.


 The combined result of applying $D$ and $ROTS$ as described above is a
 very large collection of prime alternating knots of $n+1$ crossings
 (for example, in our computation of the $40,619,385$ prime alternating
 knots of 19 crossings, $39,722,121$ were constructed by $D$ and $ROTS$).

 The next stage in our construction consists of repeatedly performing the
 following action: for each positive 2-group of each configuration in the
 collection of all knots of $n+1$ crossings that have been constructed so
 far, apply $T$ and, if the result is a new knot (which is determined
 by computing its master array and comparing it to those already in the
 collection of all knot produced so far), add it to the collection. This
 process terminates when the collection reaches a stage where subsequent
 applications of $T$ do not produce any new knots.


 Now we are ready to apply $OTS$ for the first time. We examine each knot
 configuration of each knot in the current collection to identify the $OTS$
 scenarios and to each one, we apply $OTS$. Each new knot configuration is
 added to the collection, and the enlarged collection becomes the current
 collection. When no new knot can be created by an application of $OTS$ to
 any configuration of any knot in the collection, this stage is complete.


 The construction is completed by iterating the $T$ and $OTS$ operators as
 described above, each operator only being applied to the knots produced by the
 other operator in the preceding run. Eventually, a stage is reached at
 which no new knots can be produced by applying either $T$ or $OTS$, and
 at this point, the construction of all of the prime alternating knots
 of the given crossing is complete.
 
 We now begin the process of examining the operators in the order that
 they will be applied in order to identify some of the redundant
 constructions and establish some means whereby we can prevent most of the
 unnecessary work from being done. In particular, for each operator, we
 shall establish that it suffices to apply the operator only to full group
 configurations, and that, with the possible exception of the $OTS$ operator,
 the configuration that results from an application of an operator to a full
 group configuration is again a full group configuration. Note that of the
 four operators, $OTS$ is the only one that does not operate on a group,
 and we shall see that this makes it somewhat more complicated to deal with.

\subsection{The $D$ operator}

 Observe that if $D$ is applied to two different crossings in the same
 subgroup of a negative group, the resulting configurations are identical,
 since in both cases, the size of the subgroup simply increased by 1.
 
 On the other hand, if $D$ is applied to a crossing of a positive group of
 2 or more,
 the result is to replace the crossing by a group of 2 which forms
 a min-tangle in the orbit of the group to which the crossing had belonged
 (now a group of size one smaller). It follows that the result of applying
 $D$ to two different crossings of a 
 positive group is to produce two flype equivalent configurations.

 As for the issue of full group versus split group configurations, it
 suffices to observe that if $D$ is applied to a crossing of a subgroup
 at some position of the group's orbit, and subsequently the other
 subgroups of the group are flyped to this position, the result is
 the same as first flyping the other subgroups to the selected position,
 then applying $D$ to the originally selected crossing.

 Furthermore, we observe that if $D$ is applied to a full group configuration,
 the result is still a full group configuration.
 
 Since we require a single configuration of a prime alternating knot
 from which to compute its master array, it is not necessary to apply
 $D$ to more than one crossing from a given group, nor is it necessary
 to apply $D$ to any split-group configurations. Thus we need only apply
 $D$ to full group configurations.

 Even though the decision to work only with full group configurations
 does reduce the number of configurations of a given knot that must be
 examined for potential applications of $D$, there will in general be many
 full group cofigurations of the knot. However, when we consider any
 group $G$ to which $D$ can be applied, and any flype move which does not
 move any crossings of $G$ along its orbit, then the effect of the flype
 move is restricted to some min-tangle of the orbit. It follows that
 the $D$ operation and the flype move commute with each other in such a
 case. Thus it suffices to extract the group code for a single configuration
 from the master array of the given knot, examine the selected configuration
 for a group $G$ to which $D$ is to be applied, then form all configurations
 that can be obtained by group flyping $G$ to its various orbit positions.
 Apply $D$ to the group $G$ in each of these configurations, and if the
 resulting knot is a new one, add it to the knot database.

\subsection{The $ROTS$ operator}
              
Just as for $D$, $ROTS$ operates on a group, so any flype moves that do
not involve the group will commute with $ROTS$. This observation implies
that we need only investigate full group configurations,
with the possible exception of a negative group of size at least four.
If we have a negative group of 4 or more crossings, we could flype all
but three of its crossings to a new position, thereby creating a situation
where $ROTS$ can be applied. However, upon performing $ROTS$ on the
negative 3-group, we flype the crossings back to their original
position, which results in a negative 2-group attached to the $rots$ tangle.
This same knot configuration will be obtained from the $K_A$ knot 
configuration that results from applying $D$ to the loner or negative group that is
obtained when $D^-$ is applied to the negative group that is attached to 
the $rots$ tangle. Consequently, we see that it suffices to apply $ROTS$
only to full group configurations. Furthermore, the configuration that
results from an application of $ROTS$ to a full group configuration is
itself a full group configuration.

Given a master array for a knot, it suffices to select any one full group
configuration and search for a negative group of 2 or 3. If such a
negative group is found, we shall make all possible full group
configurations from the initial one by group flyping the group to its
various positions. Then $ROTS$ 
will be applied to the group in each of these configurations. Any new 
knots that result will be added to the knot database.

\subsection{The $T$ operator}

Just as for $D$ and $ROTS$, $T$ operates on a group and so $T$ will
commute with any flype operation that does not involve any crossings
of the group that $T$ is operating on. Consider any full group
configuration that has a positive group $G$ of size at least 3.
Flype all but two of the crossings of $G$ to some common position
in the orbit of $G$, and let $G'$ denote the positive group that is
formed by these crossings in their new position.
What remains is a positive 2-group for $T$ to
operate on. In other words, by splitting this group, we have created
a chance to apply $T$ that did not exist in a full group configuration.
We shall show that the result of applying $T$ to the positive 2-group
that results when we split a positive group in this fashion will
already have been obtained by the application of either $ROTS$ (if $G$
is a 3-group) or $D$. In the
configuration that results when $T$ has been applied to the positive
2-group, the positive 2-group has been replaced by its turned image
(still a positive 2-group), which now becomes another min-tangle in the
orbit of $G'$, while $G'$ is now a negative group. Flype $G'$ back
to its original position, and denote this knot configuration by $C_1(K)$.
In $C_1(K)$, the tangle which consists of the crossing of $G'$ that is
adjacent to the turned 2-group together with the turned 2-group is a
$rots$-tangle. If $G'$ is a loner, this same configuration would have
been obtained by $ROTS$, while if $G'$ is not a loner, then it is a
negative group, and so the configuration would have been produced
by an application of $D$ to the configuration from $C_1(K)$ obtained by
reducing the size of the negative group $G'$ (that is, by applying
$D^-$ to $G'$).

It therefore suffices to apply $T$ to positive 2-groups in full group
configurations only. The result of applying $T$ to a positive
2-group $G$ in some full group configuration is a split group configuration
if and only if $G$ is a core min-tangle in the orbit of some negative
group or loner. Thus if we take care to select the full group configuration
to which $T$ is to be applied on $G$ so that if $G$ is the core tangle
of some negative group or loner, the negative group or loner is placed
in a position on one side of $G$ or the other. We shall see shortly that
rather than apply $T$ to a group code, we in fact will apply $T$ to a
master group code, and as a result, we will be in a position to instantly
detect those situations when the result of turning a particular positive
2-group will result in the creation of a larger positive group.

With this precaution, it follows that the result of
applying $T$ to a full group configuration will again be a full group
configuration.

We have also established that from each
master array, we may select just one full group configuration, and for
each positive 2-group in that configuration, we form all full group
configurations from the selected one by group flyping the positive 2-group
to each of the various positions of its orbit. We apply $T$ to the
positive 2-group in each of the these configurations, and any new knots
that are so obtained are added to the knot database.

\subsection{The $OTS$ operator}

While $OTS$ does not operate on groups, it is nevertheless
evident that $OTS$ will commute with any flype that does not involve
the three crossings in the $OTS$ 6-tangle on which $OTS$ will operate.
Furthermore, since no two crossings in an $OTS$ 6-tangle can belong to
the same group, it turns out that if $c_1$, $c_2$ and $c_3$ are the
three crossings which make up an $OTS$ 6-tangle, then performing the
$OTS$ operation on this 6-tangle will remove each of $c_1$, $c_2$ and
$c_3$ from their respective groups, but does not interfere with the
orbits of the respective groups. Thus the $OTS$ operation on this
$OTS$ 6-tangle will commute with group flypes which put the groups of
the crossings $c_1$, $c_2$ and $c_3$ into full group, with the
respective groups positioned at the locations of $c_1$, $c_2$ and
$c_3$. In other words, it suffices to only apply $OTS$ to full group
configurations. However, unlike the case for the other three operators,
from each master array we must construct each and every full group
configuration and search each one for $OTS$ 6-tangles on which to apply
$OTS$ (later on, we shall show that many of these configurations do not
have to be constructed, and even in those that do, many if not most
of the $OTS$ 6-tangles they contain will not have to be considered).

There is an interesting observation to be made at this point. The
configuration that results from an application of $OTS$ need not
be a full group configuration, which presents us with a complication, since
we wish to be working only with full group configurations. However, we
shall see shortly that one of the major benefits of master group code is
that we may apply the operator in question to a master group code, rather
than to a group code. Consequently, we shall have a simple mechanism
available to us to detect those instances when a group will become
larger as a result of an $OTS$ operation, and we shall then see that
it is extremely easy to handle this situation.

Once we have a master group code for the resulting knot, we construct
its master array and submit it to the knot database.

%% file: implem2.tex
\section{An Implementation of the Construction, Part II: Major Reductions in
the Work Done by the Operators}

As we have described the process so far, all of the operators will be 
performing a tremendous amount of redundant work. Some indication of
the redundancy is available from the run times for our calculations
of the prime alternating knots up to and including those at 19 crossings.
To facilitate the discussion, we shall refer to the two operators $D$ 
and $ROTS$ combined as $DROTS$. In the production of the 19 crossing
knots, we found that the run time for $T$ was about one half that of
$DROTS$ to complete, while the application of $OTS$ to the combined output
from $DROTS$ and $T$ took about the same amount of time as $DROTS$ and $T$
combined. As a point of interest, we note that
in every case so far (which is up to and including the knots at 19 crossings),
all knots have been produced by the end of the first application of $T$ and
$OTS$. We speculate that this will always be the case, but this is yet to be
proven.

We return to the assessment of the rate at which $T$ and $OTS$ are
producing knots. $T$ and $OTS$ account for less than 2\% of the total knot
production, yet they take over 60\% of the production time! Our objective
is to bring the amount of work done by $T$ and $OTS$ more in line
with their rate of production. In fact, it turns out that a careful
analysis of each operator permits us to dramatically reduce the amount
of work that each must do. In the production of the prime alternating
knots of 19 crossings prior to the implementation of the refinements 
to be described below, the number of knots that were built was approximately
8 times the total number of knots at the crossing size.
Preliminary calculations indicate that the refinements will bring the
ratio down to less than 2.

As we progress through these discussions, the value of the master array
will become increasingly more evident. Of course, it has tremendous value
by virtue of the fact that two prime alternating knots are flype equivalent
if and only if their master arrays are identical, but now we shall see the
returns deriving from the fact that the master array simultaneously describes
all possible (minimal crossing) configurations of the knot. The general
theme will be to establish that the four operators can be made to operate
on the master array, rather than a configuration, and this observation 
alone accounts for a monumental reduction in the amount of work to be
done in the generation of the prime alternating knots of $n+1$ crossings
from (the master arrays of) the prime alternating knots of $n$ crossings.

The following terminology will be useful for the presentation of the
reduction rules.

\begin{definition}\label{drots tangles}
 Any tangle of the form shown in Figure \ref{drots tangle}
 is called a $drots$ tangle, positive if the 2-group is positive and
 negative if the 2-group is negative.

 Moreover, in the master array of a knot, any 2-group which serves as a
 min-tangle for some group shall be called a $drots$ 2-group, with the
 same sign as the 2-group. The group in whose orbit a $drots$ 2-group
appears as a min-tangle is called the {\it orbiter} of the $drots$ 2-group.
\end{definition}

 \begin{floatingfigure}{.12\hsize}
 \centering
    $\vcenter{\xy /r27pt/:,
     (1.5,0)="a12",
     (1.5,-.5)="a13",
     (2.5,0)="a14",
     (2.5,-.5)="a15",  
     "a12";"a14"**\dir{-},
     "a13";"a15"**\crv{(1.75,-.5) & (2.5,0) & (2.25,.5) & (1.75,.5)
      & (1.5,0) & (2.25,-.5)},    
    \endxy}$
    \caption{}
  \label{drots tangle}
\end{floatingfigure}

If a 2-group is a $drots$ 2-group in a master array, then in any
configuration of the knot that has at least one of the crossings of
the orbiter group in a position adjacent to the $drots$ 2-group, the 
$drots$ 2-group together with such a crossing will form a $drots$ 
tangle. We note that a $drots$ 2-group itself has trivial orbit. 
Furthermore, any orbiter group must be either a loner or a negative group.

\begin{figure}[ht]
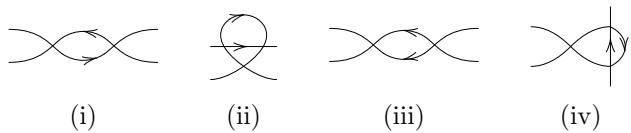

  \centering
  \begin{tabular}{c@{\hskip20pt}c@{\hskip20pt}c@{\hskip20pt}c}
   $\vcenter{\xy /r25pt/:,
   (2,0)="a1",(4.25,0)="b1","a1";"b1"**\crv{(2.5,0) & (2.75,-.5) & (3.5,-.5) & (3.75,0)}?(.6)*\dir{}="c"-<9pt,2pt>;"c"**\dir{}*\dir{>},
   (2,-.5)="a16",(4.25,-.5)="b16","a16";"b16"**\crv{(2.5,-.5) & (2.75,0) & (3.5,0) & (3.75,-.5)}?(.5)*\dir{}="c"+<9pt,-1.5pt>;"c"**\dir{}*\dir{>},
   \endxy}$ &
       $\vcenter{\xy /r25pt/:,
        (1.5,0)="a12",
        (1.5,-.5)="a13",
        (2.5,0)="a14",
        (2.5,-.5)="a15",  
        "a12";"a14"**\dir{-}?(.55)*\dir{}="c"-<9pt,.25pt>;
              "c"**\dir{}*\dir{>},
        "a13";"a15"**\crv{(1.75,-.5) & (2.5,0) & (2.25,.5) &
        (1.75,.5) & (1.5,0) & (2.25,-.5)}?(.5)*\dir{}="c"-<7pt,2pt>;
        "c"**\dir{}*\dir{>},    
       \endxy}$ &
       $\vcenter{\xy /r25pt/:,
        (2,0)="a1",(4.25,0)="b1","a1";"b1"**\crv{(2.5,0) & (2.75,-.5) &
          (3.5,-.5) & (3.75,0)}?(.5)="c"-<9pt,2pt>;
          "c"**\dir{}*\dir{<},
          (2,-.5)="a16",(4.25,-.5)="b16","a16";"b16"**\crv{(2.5,-.5) &
          (2.75,0) & (3.5,0) & (3.75,-.5)}?(.5)*\dir{}="c"+<9pt,-1.5pt>;
          "c"**\dir{}*\dir{>},
        \endxy}$ &
        $\vcenter{\xy /r30pt/:,
          (.5,.25)="q1",
          (.5,-.25)="q2",
          (1,0)="q3",
          (1.5,.5)="q4",
          (1.5,.25)="q5",
          (1.5,-.25)="q6",
          (1.5,-.5)="q7",
          "q2";"q3"**\crv{(.75,-.25)},
          "q3";"q5"**\crv{(1.25,.25)},
          "q5";"q6"**\crv{(1.75,.125) & (1.75,-.125)}?(.6)="xx",
          "xx"+(-.03,.5);"xx"**\dir{}*\dir{>},
          "q6";"q3"**\crv{(1.25,-.25)},
          "q3";"q1"**\crv{(.75,.25)},
          "q7";"q4"**\dir{-}?(.6)*\dir{>},
         \endxy}$\\
     \noalign{\vskip6pt}
        (i) & (ii) & (iii) & (iv)
   \end{tabular}
   \caption{How $drots$ tangles arise}
   \label{production of drots tangles}
 \end{figure}

Note that a positive $drots$ tangle is simply a $ROTS$ tangle, and a
$ROTS$ tangle results when the $ROTS$ operator is applied to a negative
2-group. On the other hand, a
negative $drots$ tangle results from applying $D$ to a positive 2-group.
In Figure \ref{production of drots tangles} (i) and (ii) we show a negative
2-group and the result of applying $ROTS$ to it, while in (iii) and
(iv), we show a positive 2-group and the result of applying $D$ to it.
It follows that any knot whose master array contains a $drots$ 2-group
can be obtained by an application of one of $D$ or $ROTS$.

\subsection{Reduction in the work to be done by $D$ on negative groups and loners.}

Observe that as things now stand, a knot of $n+1$ crossings with
two or more negative groups will be constructed at least once per
negative group during the application of $D$ to the various knots of
$n$ crossings, so there is considerable obvious redundancy in $D$. We
shall set out some rules which will allow $D$ to decide to pass over
an opportunity to be applied to a negative group in some knot, with
the guarantee that the knot that would result from the skipped
application of $D$ will be produced by some other application of $D$.

Suppose that we have a master array for a prime alternating knot $K$
of $n$ crossings. Let $G$ denote a negative group or a loner. We set out
the following conditions under which $D$ will be applied to $G$.
We shall let $K_1$ denote the knot in full group configuration that would
be obtained by an application of $D$ to $G$ in some full group
configuration of $K$, and let $G_1$ denote the negative group that would
be obtained from $G$ by this application of $D$.

We remark that our approach will be to work on master arrays, thereby
obtaining a master group code for the knot that results from the application
of $D$. The application of $D$ to a master array is extremely efficient, and
while we could examine the master array of the original knot to see it
meets certain conditions which determine whether or not we should apply $D$,
it is probably more efficient to actually perform $D$ and examine the
resulting master group code that results to determine whether or not to
subject it to further processing and submit it to the database of all knots
constructed so far, or else to discard it. We shall take the latter approach.

\begin{alphlist}

\item if the size of the negative group $G_1$ is greatest among all
 negative groups in $K_1$, and $G_1$ is the only negative group of this
 size in $K_1$, apply the following test.

\begin{romanlist}
 \item if the original group $G$ has a trivial flype orbit, then
  construct the master array for $K_1$ and submit it to the database.
 \item if the original group $G$ has a nontrivial flype orbit,
  then submit $K$ and $G$ to the symmetry test outlined below to determine
  in which positions of the orbit of $G$ in $K$ we must place $G$ and
  apply $D$.
\end{romanlist}

\item if the size $m$ of the negative group $G_1$ is greatest among all
 negative groups in $K_1$, but $K_1$ has more than one negative group
 of this size, let $\cal{H}$ denote the collection of all of
 the negative groups in $C(K_1)$ of size $m$ that have
 a non-trivial flype orbit, and $\cal{J}$ denote the collection of
 all of the negative groups in $C(K_1)$ of size $m$ that have a
 trivial flype orbit (except for those negative 2-groups which belong to
 a $drots$ tangle in the case $m=2$), and apply the following test.

\begin{romanlist}

\item if $G_1$ has a trivial flype orbit and $\cal{H}\ne\emptyset$, then
 do not apply $D$ to $G$ in any position of its orbit (the orbit of
 $G$ will be trivial as well if $G$ is a group of 2 or more, but if $G$
 is a loner, it is possible that $G$ could have had an orbit in the
 positive direction).
\item if $G_1$ has a trivial flype orbit and $\cal{H}=\emptyset$, but
 $\cal{J}-\{\, \hbox{$G_1$}\,\}\ne\emptyset$, then perform a group $LNR$
 competition between $G_1$ and each position of each group belonging to
 $\cal{J}-\{\, \hbox{$G_1$}\,\}$. If $G_1$ is the winner of the
 competition, then compute the master array for $K_1$ and submit it
 to the database.
\item if $G_1$ has a nontrivial flype orbit and $\cal{H}-\{\,\hbox{$G_1$}\,\}
\ne \emptyset$, then compare the size of the orbit of $G_1$ to that of
 the other groups in $\cal{H}$. If the orbit of $G_1$ is not of maximal
 size, then do not apply $D$ to $G$ in any position. On the other hand,
 if the orbit of $G_1$ is of maximal size $k$, then submit $G_1$ to a core
 $LNR$ competition (described below) against each of the groups in
 $\cal{H}- \{\,\hbox{$G_1$}\,\}$ whose orbit size is $k$ . If $G_1$ is the
 competition winner, then send $G$ off to the symmetry test described below
 to determine the positions in the orbit of $G$ to which $D$ is to be
 applied.
\end{romanlist}

In all other cases, submit $K$ and $G$ to the symmetry test described 
below, and perform $D$ on $G$ in every position indicated by the 
symmetry test, calculate the master array for the resulting knot, and 
submit it to the database.

\end{alphlist}

We shall first describe the $LNR$ competition scheme, and then present
the symmetry test.

The abbreviation $LNR$ stands for a {\it least non-repeating sequence},
by which we mean a longest possible subsequence of consecutive entries
from a specified starting point in a master group code with a specified
direction of traversal, subject to the condition that no two entries of
the sequence belong to one group (where for the purpose of this concept,
all positions of a group are considered as belonging to the group).
A group $LNR$ competition consists of the following steps. At
each end arc of each group in the competition, traverse the arc in the
direction away from the group, and examine the first group
to be encountered. Stop the competition for any strand for which the
group encountered is not the smallest (in the lexicographical ordering)
of those encountered by all of the strands for all of the groups in the
competition. If after this first stage, there is two or more strands
still active among all of the groups in the competition, and at least
one of these belongs to $G_1$, carry out the next stage, whereby
each strand that is still in the competition is traversed to the next
group to be encountered, and the assessment is repeated (but now, a strand
will also be stopped if the group encountered is a repetition of one
encountered earlier on the strand, since that would cause the sequence of
groups encountered along the strand to fail to be a $LNR$). The process
stops when $G_1$ no longer has any strands in the competition, but there
are still active strands, in which case $G_1$ is declared to have lost the
competition, or until all strands $G_1$ have stopped, and $G_1$ had at
least one strand active up to the end, in which case $G_1$ has won (or
tied, but we declare a tie to be a win) the competition.

A core $LNR$ competition is a minor variant of the group $LNR$ competition
described above. The competition is applicable only to negative groups with
a non-trivial flype orbit, and instead of starting at the group ends and 
following the strands out from the group, we go to the core tangle for the 
group, and follow the outbound strands from the core tangle through the copy
of the group that sits on either side of the core tangle, taking the
resulting location as our starting point for each strand in the competition.
The direction of travel on the strand will be away from the core tangle.

We now describe the symmetry test. The underlying observation is that
a group $G$ may have two or more positions, but it may happen that
as a result of symmetry within one or more of the min-tangles in the
orbit of $G$, if $D$ is applied to $G$ in two different positions in
its orbit, the resulting knots are flype equivalent. Accordingly, we
now proceed to identify certain such situations.

\begin{definition}\label{rotating a min-tangle}
Let $C(K)$ be a full group configuration of a prime alternating knot $K$,
$G$ be a group in $C(K)$, and let $T$ be a min-tangle
in the orbit of $G$. The following process shall be referred to as
{\em rotating} $T$.  Cut the four arcs incident to $T$ and label each cut
arc $i$ with two labels, one on either side of the cut so that the side
that is no longer incident to $T$ is labelled by $i'$ and the side that is still
incident to $T$ is labelled by $i''$. Then rotate $T$ as if it were
involved in a flype of $G$ from the pair of arcs on one side of $T$
to the pair of arcs on the other side of $T$, and connect the cut
arcs in such a way that if the arcs on a side were labelled $i$ and $j$, 
then $i'$ is connected to $j''$ and $j'$ is connected to $i''$. 
\end{definition}

Since one of the two arcs on one side of $T$ represents an overpass
of $G$, while the other arc on the same side represents an underpass of
$G$, it follows that rotating a min-tangle of a group in a configuration
of a prime alternating knot results in a configuration of an alternating
knot, necessarily prime.

\begin{definition}
Let $C(K)$ be a full group configuration of a prime alternating knot $K$,
$G$ be a group in $C(K)$, and $T$ be a min-tangle
in the orbit of $G$. Let $K_1$ denote the prime alternating knot whose
configuration is obtained from $C(K)$ by rotating $T$. $T$ is said to
be {\em flype-symmetric} if and only if $K_1=K$.
\end{definition}

The significance of a flype-symmetric min-tangle $T$ in the orbit of $G$,
where $G$ is either a negative group or else a loner with an
orbit in the negative direction, is that applying $D$ to $G$ whether $G$
is adjacent to $T$ on one side or the other results in the same knot.

\begin{definition}
 Let $K$ be a prime alternating knot with master array $M$, and let
 $G$ be a group in $C(K)$. Any choice of two positions in the orbit of $G$
 provides a splitting of the orbit of $G$ into two tangles, each referred to
 as the complement of the other (if $T$ denotes one of the two tangles, we
 shall denote the other one, the complement of $T$, by $T^c$). Suppose
 that two positions in the orbit of $G$ have been chosen. Let $T$ denote
 one of the two tangles that results. $M$ determines an orientation of $K$,
 and according to this orientation, two of the four arcs incident to $T$
 are entering $T$, while the other two are exiting $T$. Let $e_{in}$ denote
 an arc entering $T$ and let $e_{out}$ denote the arc by which $T$ is first
 exited after $T$ has been entered by $e_{in}$. The sequence $s$ in $M$
 that is delimited by $e_{in}$ and $e_{out}$ is called a {\em master array
 sequence determined by $T$}. Let $[s]$ denote the sequence that is obtained
 from $s$ by traversing $s$, relabelling the groups so that the $i^{th}$ group
 of size $k$ that is encountered in the traversal of $s$ receives label $k_i$
 (with the correct sign), and the positions of group $k_i$ are labelled in the
 order that they are encountered in the traversal of $s$. Furthermore, let
 $[s]^r$ denote the sequence that is obtained by applying the same procedure
 to the same sequence of the master array, but traversed in the opposite
 direction.
\end{definition}

\begin{proposition}\label{identify flype-symmetry}
 Let $C(K)$ be a configuration of a prime alternating knot $K$, $G$ a 
 group in $C(K)$, and $T$ a min-tangle in the orbit of $G$. Further let
 $s$ and $t$ denote the two master array sequences that are determined
 by $T$. Then $T$ is flype-symmetric if:

 \begin{alphlist}
  \item $G$ is negative, $T$ is the core of $G$'s orbit, $[s]=[s]^r$, and
    $[t]=[t]^r$, or
  \item $G$ is negative, $T$ is not the core of $G$'s orbit, and $[s]=[t]^r$,
    or
  \item $G$ is positive, and $[s]=[t]$.
 \end{alphlist}
\end{proposition}

\begin{proof}
  Suppose $G$ is negative, $T$ is the core tangle of $G$'s orbit, and
  $[s]=[s]^r$, $[t]=[t]^r$. Then upon rotating $T$, the segments $s$ and
  $t$ in the master array $M$ get replaced by $s$ reversed and $t$ reversed,
  respectively to produce a master group code of the resulting knot. Since
  $[s]=[s]^r$ and $[t]=[t]^r$, the master array that results from this
  master group code is again $M$, whence $T$ is flype-symmetric.

  Now suppose that $G$ is negative, $T$ is not the core tangle of $G$'s
  orbit, and $[s]=[t]^r$. Upon rotating $T$, the segments $s$ and $t$ in
  $M$ are replaced by $t$ reversed and $s$ reversed, respectively,
  providing a master group code for the resulting knot. Since $[s]=[t]^r$,
  we also have $[t]=[s]^r$, and so the master array that is obtained from
  this master group is $M$, whence $T$ is flype-symmetric.

  Finally, suppose that $G$ is positive, and that $[s]=[t]$. Upon rotating
  $T$, the segments $s$ and $t$ in $M$ get replaced by $t$ and $s$,
  respectively, providing a master group code for the resulting knot. Since
  $[s]=[t]$, the master array that is obtained from this master group code
  is again $M$, so $T$ is flype-symmetric.
\end{proof}


  The way in which the above information is utilized in the application of
  $D$ on a negative group $G$ or a loner with an orbit in the negative
  direction is as follows. After having decided that we shall apply $D$
  to the negative group $G$ in the knot $K$, we locate $G$ in its zero
  position in the master array of $K$, and mark this position as one in
  which we shall perform the $D$ operation. We then traverse the knot
  from this position in the orbit of $G$ in the direction determined by
  the master array and examine the first min-tangle $T$ that is
  encountered. If $T$ is flype-symmetric, then we do not need to place $G$
  on the arcs on the other side of $T$ and perform the $D$ operation on $G$
  in this position. We check the conditions outlined in \ref{identify
  flype-symmetry} and if met, we know that $T$ is flype-symmetric and we
  move on to repeat the process with the next min-tangle in the orbit of
  $G$. If the conditions are not met, we mark this location as one where
  we shall perform the $D$ operation on $G$. We then continue on to
  investigate the next min-tangle. This process is repeated until the orbit
  has been completely examined. For each position in the orbit of $G$ that
  has been marked by the procedure outlined above, $D$ is applied to $G$ in
  that position, and the master array of the resulting knot is constructed
  and the knot is submitted to the database.

\subsection {Reduction in the work to be done by $D$ on positive 2-groups}
  If a prime alternating knot $K$ contains a positive 2-group $G$ such that
  if $D$ was to be applied to $G$, the resulting knot $K_1$ would contain
  either a negative 2-group that is not in a negative $drots$ tangle, or a
  negative group of size at least three, then the resulting knot would
  already have been made by $D$ on negative groups, and so we shall not
  apply $D$ to such a positive 2-group. We remark that it suffices to
  check $G$ in any position of its orbit in order to determine whether or
  not any of these conditions are met. On the other hand, suppose that 
  $K_1$ contains no negative groups other than negative 2-groups in negative
  $drots$ tangles. We may still opt to defer to one of these by setting
  up a $LNR$ competition between $G_1$ and the other negative 2-groups.
  As before, if $G_1$ wins the competition, then we construct the master
  array for $K_1$ and submit it to the database, otherwise we choose not
  to process the application of $D$ to $G$.

\subsection{Reductions in the work to be done by $ROTS$}

\noindent {\bf First Reduction for $ROTS$.} If a master array contains
two or more negative groups, then the master array will not be
submitted to $ROTS$. For the result of applying $ROTS$ to a negative
2-group or a negative 3-group would still leave a negative group 
elsewhere in the knot, which means that the resulting knot would
already have been produced by $D$. 

\noindent {\bf Second Reduction for $ROTS$.} If a master array contains
just one negative group, and that one negative group is the negative
2-group of a negative $drots$-tangle, shown in Figure \ref{redundant
rots} (i), then this master array will not be submitted to $ROTS$.
For if we were to apply $ROTS$ to this
negative 2-group, the result would be as shown in Figure \ref{redundant
rots} (ii). But this last tangle is exactly what
would result from an application of $D$ to $c$ (where $c$ is either a
loner, or a crossing in a negative group). It follows therefore that
if we were to apply $ROTS$ to such a negative 2-group, the knot that
would result would already have been constructed by $D$, and so this
would be redundant work. Thus we shall never apply $ROTS$ in this
situation.

 \begin{figure}[ht]
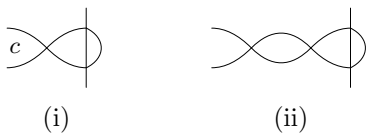

  \centering
  \begin{tabular}{c@{\hskip40pt}c}
   $\vcenter{\xy /r30pt/:,
       (.5,.25)="q1",
       (.5,-.25)="q2",
       (1,0)="q3",
       (1.5,.5)="q4",
       (1.5,.25)="q5",
       (1.5,-.25)="q6",
       (1.5,-.5)="q7",
       (.6,0)*!<0pt,0pt>{\hbox{\small $c$}},
       "q2";"q3"**\crv{(.75,-.25)},
       "q3";"q5"**\crv{(1.25,.25)},
       "q5";"q6"**\crv{(1.75,.125) & (1.75,-.125)},
       "q6";"q3"**\crv{(1.25,-.25)},
       "q3";"q1"**\crv{(.75,.25)},
       "q4";"q5"**\dir{-},
       "q5";"q6"**\dir{-},
       "q6";"q7"**\dir{-},
      \endxy}$
      &
      $\vcenter{\xy /r30pt/:,
        (4,.25)="q8",
        (4,-.25)="q9",
        (4.5,0)="q10",
        (5.25,0)="q11",
        (5.75,.5)="q12",
        (5.75,.25)="q13",
        (5.75,-.25)="q14",
        (5.75,-.5)="q15",
        "q8";"q10"**\crv{(4.25,.25)},
        "q10";"q11"**\crv{(4.75,-.25) & (5,-.25)},
        "q11";"q13"**\crv{(5.5,.25)},
        "q13";"q14"**\crv{(6,.125) & (6,-.125)},
        "q14";"q11"**\crv{(5.5,-.25)},
        "q11";"q10"**\crv{(5,.25) & (4.75,.25)},
        "q10";"q9"**\crv{(4.25,-.25)},
        "q12";"q13"**\dir{-},
        "q13";"q14"**\dir{-},
        "q14";"q15"**\dir{-},
       \endxy}$\\
       \noalign{\vskip6pt}
       (i) & (ii) 
      \end{tabular}
      \caption{A second redundant $ROTS$ scenario}
    \label{redundant rots}
  \end{figure}

\subsection{Reductions in the work to be done by $T$ and $OTS$}
Our first step is to remove from the $DROTS$ output those master arrays
that we can guarantee will not produce any new knots by an application of
either $T$ or $OTS$.

Our first reduction stems from the observation that if a master array
that was produced by $DROTS$ contains a negative $drots$ 2-group, then
since the orbiter group of the $drots$ 2-group is either
negative or else a loner, it follows that any application of $T$
will not operate on the $drots$ 2-group nor on its orbiter group, nor
will the orbit of the orbiter group be affected. Thus the master array of
any knot that is produced by an application of $T$ to any configuration
that is obtained from a master array that contains a negative $drots$
2-group will contain a $drots$ 2-group (which might be positive since
the application of $T$ could change the sign of one or more groups in the
configuration). Thus all knots that $T$ can produce from such a master
array will already have been produced by $DROTS$. 

Next, consider the result of applying $OTS$ to a configuration
obtained from a master array which contains a negative $drots$ 2-group.
If $OTS$ is applied to an $OTS$ 6-tangle which does not involve either
of the crossings of the $drots$ 2-group, then the resulting configuration
has a negative 2-group. In such a case, the knot that the 
resulting configuration represents will already have been produced by $D$.
Suppose now that there is an $OTS$ 6-tangle which involves some of the
crossings of the $drots$ 2-group. Since no $OTS$ 6-tangle can contain two
crossings from the same group, this $OTS$ 6-tangle must contain exactly
one crossing from the negative $drots$ 2-group. The following proposition
will allow us to deal with this situation.

\begin{proposition}\label{orbiter and ots}
 Let $C(K)$ be a full group configuration of a prime alternating knot $K$
 and suppose that $G_D$ is a $drots$ 2-group (positive or negative) such
 that there is an $OTS$ 6-tangle that contains a crossing of $G_D$. Then
 the orbiter group $G$ of $G_D$ must be a loner adjacent to $G_D$.
\end{proposition}

\begin{proof}
  Let $c$, $x$ and $y$ be crossings which form an $OTS$ 6-tangle, and
  suppose that $c$ is in $G_D$ (so neither $x$ nor $y$ can belong to
  $G_D$). Suppose that neither $x$ nor $y$ is a crossing of $G$. Since
  $G_D$ is a min-tangle in the orbit of $G$, $x$ and $y$ belong to
  min-tangles $T_1$ and $T_2$, respectively, in the orbit of $G$. 
  Now $T_1$ and $T_2$ are adjacent to $G_D$, one on each side, so the
  group $G$ lies between $T_1$ and $T_2$. Thus there can be no edge directly
  joining $x$ to $y$, so this situation can't arise. Accordingly, at 
  least one, and therefore exactly one of $x$ and $y$ belong to $G$,
  whence $G$ is adjacent to $G_D$.  Suppose that $x$ is in $G$. Then
  $y$ is in a min-tangle in the orbit of $G$, on the other side of $G_D$
  from $G$. Thus for $x$ and $y$ to be adjacent, it must be that $G$
  consists only of $x$; that is, $G$ is a loner adjacent to $G_D$, as
  required.  
\end{proof}

\begin{figure}[ht]
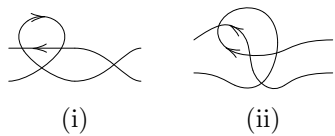

  \centering
  \begin{tabular}{c@{\hskip20pt}c}
   $\vcenter{\xy /r25pt/:,
      (1.5,0)="a12",
      (1.5,-.5)="a13",
      (2.5,0)="a14",
      (2.5,-.5)="a15",  
      "a12";"a14"**\dir{-}?(.55)*\dir{}="c"-<9pt,.25pt>;"c"**\dir{}?(.5)*\dir{<},
      "a13";"a15"**\crv{(1.75,-.5) & (2.5,0) & (2.25,.5) & (1.75,.5) &
        (1.5,0) & (2.25,-.5)}?(.48)*\dir{}="c"-<7pt,.45pt>;"c"**\dir{}*\dir{>},  
      (2.5,0)="a2",(3.5,-.5)="b2","a2";"b2"**\crv{(3,0) & (3.25,-.5)},
      (2.5,-.5)="a3",(3.5,0)="b3","a3";"b3"**\crv{(3,-.5) & (3.25,0)},
     \endxy}$ &
   $\vcenter{\xy /r25pt/:,
      (8,-.25)="a14",
      (.9,0)="a12",
      (.9,-.5)="a13",
      (3,0)="a14",
      (3,-.5)="a15",  
      "a12";"a15"**\crv{(1.25,.25) & (1.5,.25) & (1.75,0) & (1.75,-.5) &
          (2,-.75) & (2.25,-.75) & (2.5,-.5)}?(.15)*\dir{}="c"-<9pt,-6pt>;
          "c"**\dir{}*\dir{<},
      "a13";"a14"**\crv{(1.25,-.5) & (1.5,-.75) & (2,-.75) & (2.25,0) &
      (2,.5) & (1.5,.5) & (1.25,.25) & (1.25,0) & (1.5,-.25) & (2.25,-.25) &
      (2.5,0)}?(.75)*\dir{}="c"+<7pt,-3.5pt>;"c"**\dir{}*\dir{>},   
     \endxy}$ \\
   \noalign{\vskip6pt}
   (i) & (ii) 
  \end{tabular}
  \caption{Negative $drots$ tangles and $OTS$}
  \label{ots and negative drots tangles}
\end{figure}

 Thus in our current situation, the orbiter group must be a loner in a
 position adjacent to the $drots$ 2-group, and the $OTS$ 6-tangle is as
 shown in Figure \ref{ots and negative drots tangles} (i). The result
 after $OTS$ has been performed is shown in (ii). Note that there is a
 $drots$ tangle in (ii), so the resulting knot would already have been
 built by $DROTS$. 

 Thus any master array which contains a negative $drots$ tangle will not
 produce any new knots from applications of either $T$ or $OTS$, and so we
 shall remove such master arrays from the collection that will be given
 to $T$ and subsequently to $OTS$.

 So now let us examine the situation when we apply the $T$ operator to the
 master arrays that have been produced by $DROTS$ and which are without
 negative $drots$ 2-groups.

 We examine each such master array in turn. 

\vskip6pt
\noindent {\bf First Reduction for $T$.} Consider the case of a master array
 which contains at least two positive $drots$ 2-groups. The master array
 of any knot that results from an application of $T$ to any full group
 configuration of such a knot will necessarily contain a $drots$ 2-group
 (possibly negative, since applying $T$ can cause changes in the sign of
 a group) and so it would already have been produced by $DROTS$. Thus we
 may skip any such master array.

\vskip6pt
\noindent {\bf Second Reduction for $T$.} It is a consequence of
 Proposition \ref{kb to kc not by t}\immediate\write16{sar (implem2.tex):
 confirm this reference number} of [1] that it is not necessary
 to iterate $T$ in the application of $T$ immediately after $DROTS$.

\vskip6pt
\noindent {\bf Configuration dependent reductions for $T$.}
 Any master array that does not get set aside as a result of the 
 reduction described above is now examined for configuration dependent
 reductions. If the master array contains exactly one positive $drots$ 
 2-group, then even though it may contain other positive 2-groups, the 
 only positive 2-group which might lead to a new knot upon application 
 of $T$ is the positive $drots$ 2-group (which necessarily has trivial 
 orbit). Even this application may not be necessary, and we shall
 submit such a knot to the following test before applying $T$--the point
 is that there is only one group to examine in this particular situation). 

 If the master array contains at most one positive $drots$ 2-group, then
 for each positive 2-group $G$, unless there is a positive $drots$ 2-group,
 in which case we restrict $G$ to just the positive $drots$ 2-group,  we
 check to see if at least one of the following conditions holds.

\begin{alphlist}
 \item there is a negative group both of whose arcs are
  in the same section of $G$.
 \item there is a positive group with an arc in each section of $G$.
 \item turning $G$ results in a positive group $G'$ having the turned
  $G$ as a subgroup, but there is a positive group $H$ of size larger than
  that of $G'$ in the resulting knot.
 \end{alphlist}
 
 If any of these conditions hold for the group $G$ under consideration,
 then we do not apply $T$ to $G$. For if either of the first two
 conditions holds, then the configuration that would result from turning
 $G$ would have a negative group, so it would already have been produced
 by $D$. On the other hand, suppose that the last condition holds. Then
 the configuration $C(K)$ that results from turning $G$ would also be
 produced by an application of $T$ to the knot $C(K_1)$ obtained from $C(K)$
 by turning an end 2-subgroup of the positive group $H$ (which must be of
 size at least 3) to create a positive 2-group $H'$ in a $ROTS$ tangle.
 Thus $C(K_1)$ would be produced by $DROTS$, and $C(K)$ would be produced
 from $C(K_1)$ by turning $H'$, producing a unique positive group of
 largest size. Our convention requires that we not apply $T$ to $G$, in
 deference to this other application of $T$.

 For each positive 2-group $G$ which is to have $T$ applied to it,
 we apply $T$ to the group in each of its positions, and for each of the
 resulting knots, we construct the master array and submit it to the
 database.

 When the first round of $T$ applications is complete, we take the
 master arrays that were produced by $DROTS$ but which have no negative
 $drots$ 2-groups, together with any new knots that were produced by $T$
 as the set of knots that $OTS$ is to work on. We remark that any knots
 that $T$ produces will not have any negative (or positive) $drots$
 2-groups, so it is not necessary to examine them for the purpose of
 eliminating any that might have a negative $drots$ 2-group.

 However, before we submit any knots to $OTS$, we perform a further round
 of reductions that are unique to $OTS$.

 \begin{floatingfigure}{.2\hsize}
  \begin{tabular}{c@{\hskip20pt}c}
   $\vcenter{\xy /r40pt/:,
      (.5,0)="a1",
      (.75,0)="a2",
      (.125,1.125)="a3",
      (1.125,1.125)="a4",  
      "a3";"a4"**\dir{-}?(.525)*\dir{}="c"-<9pt,.25pt>;"c"**\dir{}*\dir{>},
      "a1";"a2"**\crv{(.5,.25) & (.75,.5) & (.75,.75) & 
        (.375,1.125) & (.375,1.375) & (.625,1.5) & (.875,1.375) & 
        (.875,1.125) & (.5,.75) & (.5,.5) & (.75,.25)}
        ?(.5)*\dir{}="c"-<7pt,2pt>;"c"**\dir{}*\dir{>}, 
    \endxy}$
  \end{tabular}
  \caption{}            
  \label{positive drots and negative group}
\end{floatingfigure}

\vskip6pt
\noindent {\bf First $OTS$ reduction:} any master array that contains a
 positive $drots$ 2-group whose orbiter group is of size at least 2 can be
 skipped by $OTS$. To see why, let $G_D$ denote a positive 2-group whose
 orbiter group $G$ has size at least 2 (see Figure \ref{positive drots 
 and negative group}).

 Observe first of all that any $OTS$ 6-tangle that does not involve
 a crossing from $G_D$ can involve at most one crossing from $G$, and the
 result of applying $OTS$ to such an $OTS$ 6-tangle would be a configuration
 which (is flype equivalent to a configuration which) has a $ROTS$ tangle
 and would therefore have been produced by $DROTS$. On the other hand, if
 there is an $OTS$ 6-tangle which involves a crossing of $G_D$, then
 by Proposition \ref{orbiter and ots}, the orbiter would be a loner.
 Since this is not the case, there can be no $OTS$ 6-tangle which involves
 a crossing of $G_D$.

So we are considering only master arrays in which any positive $drots$
2-group has a loner as its orbiter group.

\vskip 6pt
\noindent {\bf Second $OTS$ reduction:} if a master array contains a negative
group of size at least 3, a positive group of size at least 4, or more 
than two positive $drots$ 2-groups, or
if the sum of the number of negative 2-groups, positive 3-groups and 
positive $drots$ 2-groups exceeds 3, then no $OTS$ operation can destroy
all of the scenarios which would allow the resulting knot to be obtained
from $DROTS$. For if an $OTS$ 6-tangle does not involve the crossing of
the loner orbiter of a positive $drots$ 2-group (in whatever configuration),
or does not involve a crossing of a negative group of size at least 2, or
a crossing of a positive group of size at least 3, then after $OTS$
operates on this $OTS$ 6-tangle, there will still be a positive $drots$
tangle (possibly after flyping), a negative group of size at least 2, or
a positive group of size at least 3, respectively, in the configuration
that results. Such a configuration could be obtained by $ROTS$ in the case
of a positive $drots$ tangle, $D$ in the case of a negative group, or
$T$ after $ROTS$ on a negative 2-group, in the case of a positive group of
size at least 3. In any event, such a knot would already have been made
before $OTS$ begins, so any such master array can be skipped by $OTS$.

At this point, we have eliminated all master arrays which contain
more than two positive $drots$ 2-groups.

\begin{figure}[ht]
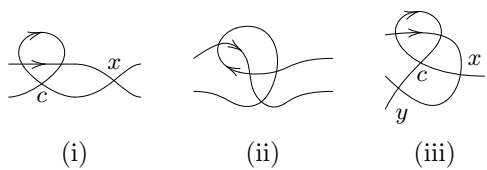

  \centering
  \begin{tabular}{c@{\hskip20pt}c@{\hskip20pt}c}
   $\vcenter{\xy /r25pt/:,
    (1.5,0)="a12",
    (1.5,-.5)="a13",
    (2.5,0)="a14",
    (2,-.5)*{\hbox{\small $c$}}, 
    (3.1,0)*{\hbox{\small $x$}},
    (2.5,-.5)="a15",  
    "a12";"a14"**\dir{-}?(.55)*\dir{}="c"-<9pt,.25pt>;"c"**\dir{}*\dir{>},
    "a13";"a15"**\crv{(1.75,-.5) & (2.5,0) & (2.25,.5) & (1.75,.5) & 
     (1.5,0) & (2.25,-.5)}?(.5)*\dir{}="c"-<7pt,2pt>;"c"**\dir{}*\dir{>},  
    (2.5,0)="a2",(3.5,-.5)="b2","a2";"b2"**\crv{(3,0) & (3.25,-.5)},
    (2.5,-.5)="a3",(3.5,0)="b3","a3";"b3"**\crv{(3,-.5) & (3.25,0)},
   \endxy}$
   &
   $\vcenter{\xy /r25pt/:,
     (8,-.25)="a14",
     (.9,0)="a12",
     (.9,-.5)="a13",
     (3,0)="a14",
     (3,-.5)="a15",  
     "a12";"a15"**\crv{(1.25,.25) & (1.5,.25) & (1.75,0) & (1.75,-.5) 
        & (2,-.75) & (2.25,-.75) & (2.5,-.5)}?(.25)*\dir{}="c"-<9pt,-6pt>;
        "c"**\dir{}*\dir{>},
     "a13";"a14"**\crv{(1.25,-.5) & (1.5,-.75) & (2,-.75) & (2.25,0) 
      & (2,.5) & (1.5,.5) & (1.25,.25) & (1.25,0) & (1.5,-.25) 
      & (2.25,-.25) & (2.5,0)}?(.75)*\dir{}="c"+<7pt,-3.5pt>;
      "c"**\dir{}*\dir{>},   
     \endxy}$ 
     & 
     $\vcenter{\xy /r25pt/:,
       (1.5,.125)="a12",
       (1.5,-1)="a13",
       (1.5,-.5)="a14",
       (3,-.5)="a15",
      (2.05,-.5)*{\hbox{\small $c$}}, 
      (2.85,-.25)*{\hbox{\small $x$}},
      (1.75,-1.1)*{\hbox{\small $y$}},
       "a12";"a14"**\crv{(2.5,.25) & (2.75,-.25) & (2.5,-1) & (2,-1)}?(.075)*\dir{}="c"-<9pt,.25pt>;"c"**\dir{}*\dir{>},
       "a13";"a15"**\crv{(1.75,-.5) & (2.5,0) & (2.25,.5) & (1.75,.5) & (1.5,0) & (2.25,-.5)}?(.5)*\dir{}="c"-<7pt,2pt>;"c"**\dir{}*\dir{>},    
       \endxy
%
     }$       
\\
     \noalign{\vskip6pt}
   (i) & (ii) & (iii) 
  \end{tabular}
  \caption{An $OTS$ scenario involving a positive $drots$ 2-group}
\label{ots and drots}
\end{figure}

\vskip 6pt
\noindent {\bf Third $OTS$ reduction:}
Consider a master array which contains a positive $drots$ 2-group in
the setting illustrated in Figure \ref{ots and drots} (i). This is
intended to show that the loner $c$ is the orbiter group of the 
positive $drots$ 2-group, and the crossing $x$ (which is either a loner
or a crossing in a negative group of 2) has a position
in its orbit which has it adjacent to the positive $drots$ 2-group
as shown in (i).

Then $c$ has two positions in its orbit, the one illustrated in (i)
and its position after flyping across the $drots$ 2-group. Both positions
of $c$ lead to the same tangle as shown in (i). Any $OTS$ 6-tangle 
must involve $c$, in whatever position $c$ is in, since otherwise the
$OTS$ operation will result in a configuration that either has a
$drots$ tangle or is flype equivalent to a configuration that has a
$drots$ tangle, and the knot represented by such a configuration has already
been produced by $DROTS$. If $x$ is not adjacent to the $drots$ 2-tangle,
as shown in (ii) (or on the other side of the $drots$ 2-group),
then $c$ is not in any $OTS$ 6-tangle and so we do not need to perform
any $OTS$ operation on such a configuration. If $x$ is in the position
shown in (i), there is the obvious $OTS$ 6-tangle that involves $c$, $x$
and one of the crossings of the $drots$ 2-group. The result of performing
the $OTS$ operation on this $OTS$ 6-tangle is shown in (ii). The resulting
configuration contains a negative 2-group, and thus the knot it represents
has already been produced by $DROTS$. 
There is one other possible $OTS$ 6-tangle that could be formed with $c$.
If the positive $drots$ 2-group is actually in the setting shown in (iii),
then $x$ has only two positions in its orbit, the position shown in (iii)
and its position after it is flyped across the $drots$ tangle (every
configuration of this knot will actually have $c$ adjacent to the positive
$drots$ 2-group to form a positive $drots$ tangle). If the $OTS$ 
operation is performed on the $OTS$ 6-tangle $c$, $x$ and $y$ as shown
in (iii), the resulting configuration has a positive 3-group, hence is
also obtainable by $T$ following $ROTS$ on a negative 2-group. On the other
hand, if we have a configuration with $y$ still adjacent to $x$, but
flyped across the $drots$ 2-group, then $y$ is no longer adjacent to 
$c$. Since we need only perform an $OTS$ operation on an $OTS$ 6-tangle
which contains $c$, there is nothing to do for this configuration. If in
addition, we flype $c$ across the positive $drots$ 2-group, then again
we find two $OTS$ 6-tangles involving $c$, one of which results in a
negative 2-group after the $OTS$ operation, the other resulting in a
positive 3-group after the $OTS$ operation. Thus in every case, the
resulting knot will already have been produced by $DROTS$ and possibly $T$.
It follows that it is not necessary to apply $OTS$ to any configuration
obtained from a master array which contains a positive $drots$ 2-group
in the setting shown in (i). Since the positive $drots$ 2-groups of this
type are so important for their role in the reduction of the work done
by $OTS$, we shall give them a name as well.

\begin{definition}\label{tight drots 2-group}
 A positive $drots$ 2-group whose orbiter group is a loner 
 and for which the $drots$ 2-group together with its
 orbiter form a min-tangle in the orbit of some group $G$ is called a
 {\it tight $drots$ 2-group}. A tangle which consists of the tight
 $drots$ 2-group, its loner orbiter and one crossing from $G$ (as shown
 in Figure \ref{ots and drots} (i)\,) is called a {\it tight $drots$
 tangle}.
\end{definition} 

Up to this point, we have determined that we need not submit any master array
from $DROTS$ to $OTS$ if it contains any of the following:

\dimen20=\baselineskip
\vskip6pt\hskip15pt
\begin{tabular}{rl}
(i) & a negative $drots$ 2-group,\\
(ii) & a positive $drots$ 2-group whose orbiter group is not a loner,\\
(iii) & a negative group of size at least 3,\\
(iv) & a positive group of size at least 4,\\
(v) & 3 or more positive $drots$ 2-groups,\\
(vi) & \vtop{\baselineskip=\dimen20\dimen0=\hsize\advance\dimen0 by -\leftskip\advance\dimen0 by -70pt
\hsize=\dimen0\noindent any combination of four or more of negative 2-groups, positive 3-groups,
 and positive $drots$ 2-groups,}\\
 (vii) & a tight $drots$ 2-group.
\end{tabular} 
\vskip6pt 

\noindent {\bf Fourth $OTS$ reduction:} At this point, the master arrays that
we are considering have at most two positive $drots$ 2-groups. We examine
those which have exactly two positive $drots$ 2-groups (each necessarily with
orbiter a loner). The only $OTS$ 6-tangles that we need consider must involve
the two orbiter loners (each in either position of its orbit). Thus if no
configuration of the knot has the two orbiter loners adjacent to each other,
there will be no $OTS$ 6-tangle to consider, and so we may forego the
application of $OTS$ to any configuration obtained from such a master array.

Suppose now that we have a master array with exactly two positive $drots$
2-groups such that the two orbiter loners are adjacent in some configuration.
We consider the various possible scenarios.

 \begin{figure}[ht]
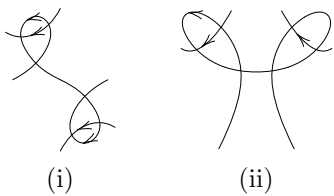

  \centering
  \begin{tabular}{c@{\hskip25pt}c@{\hskip25pt}c}
    $\vcenter{\xy /r9pt/:,
      (-2,0);(2,0)**\crv{ (-1,.5) & (0,2) & (-1,3) & (-2,2) & (-1,.5) & (0,0)
       & (1,-.5) & (2,-2) & (1,-3) & (0,-2) & (1,-.5) }?(.24)="x"?(.8)="x1",
       "x"+(.04,.01);"x"**\crv{}*\dir{>},
       "x1"+(.04,.02);"x1"**\crv{}*\dir{>},
      (-2.3,2);(0,3)**\crv{ (-1.5,1.5) & (-.5,2) }?(.42)="y",
      "y"+(.04,.015);"y"**\crv{}*\dir{>},
      (2.3,-2);(0,-3)**\crv{ (1.5,-1.5) & (.5,-2) }?(.62)="y1",
      "y1"+(.04,.015);"y1"**\crv{}*\dir{>},
    \endxy}$
    &
    $\vcenter{\xy /r38pt/:,
     (1.25,0)="a1",
     (.5,0)="a2",
     (.125,.875)="a3",
     (.625,1.375)="a4",  
     (1.125,1.375)="a5",
     (1.625,.875)="a6",
      (.125,1);(.625,1.375)**\crv{ (.275,.9) & (.475,1.1) }?(.42)="y",
      "y"+(.04,.035);"y"**\crv{}*\dir{>},
      (1.125,1.375);(1.625,1)**\crv{ (1.275,1.1) & (1.475,.9) }?(.32)="y",
      "y"+(.04,-.04);"y"**\crv{}*\dir{>},      
     "a1";"a2"**\crv{(1,.5) & (1,1) & (1.375,1.375) & (1.625,1.375) & (1.625,1.125) & (1.25,.75) & (.5,.75) & (.125,1.125) & (.125,1.375) & (.375,1.375) & (.75,1) & (.75,.5)}?(.75)*\dir{}="c"+<3pt,-.5pt>;"c"**\dir{}*\dir{>}, 
    \endxy}$
  \\
   \noalign{\vskip5pt}
  (i) & (ii)    
  \end{tabular}
  \caption{\vtop{\dimen0=\hsize\advance\dimen0 by -1.3in\hsize=\dimen0\leftskip=0pt
   \noindent The two possible $OTS$ scenarios involving two positive $drots$ 2-groups
    with adjacent orbiters}}
\label{2 drots a}
\end{figure}

First of all, suppose that the master array contains exactly two positive
$drots$ 2-groups such that some configuration of the knot has the tangle shown
in Figure \ref{2 drots a} (i). By virtue of the earlier reduction rules, we
know that the knot does not contain any tight $drots$ 2-groups, so we can be
assured that the tangle is not as shown in either Figure \ref{2 drots b} (i)
or (ii). Thus the two orbiter loners do not participate in an $OTS$ 6-tangle
in any configuration where they are adjacent to each other and each is
adjacent to its $drots$ 2-group. We must consider the possibility that in
some other configuration, the two orbiter loners might participate in an
$OTS$ 6-tangle; that is, it might be possible to start with the two positive
$drots$ 2-groups with their loners in the position as shown in Figure
\ref{2 drots b} (i) and flype one or both of the orbiters to create an
$OTS$ 6-tangle which contains the two orbiters. 

 \begin{figure}[ht]
  \centering
  \begin{tabular}{c@{\hskip25pt}c}
    $\vcenter{\xy /r9pt/:,
      (-2,0)="start",
      (2,0)="stop",
      "start";"stop"**\crv{ (-1,.5) & (0,2) & (-1,3) & (-2,2) & (-1,.5) & (0,0)
       & (1,-.5) & (2,-2) & (1,-3) & (0,-2) & (1,-.5) }?(.24)="x"?(.8)="x1",
       "x"+(.04,.01);"x"**\crv{}*\dir{>},
       "x1"+(.04,.02);"x1"**\crv{}*\dir{>},
       (-2.3,2)="TopArcStart",
       (0,3)="TopArcStop",
      "TopArcStart";"TopArcStop"**\crv{ (-1.5,1.5) & (-.5,2) }?(.42)="y",
      "y"+(.04,.015);"y"**\crv{}*\dir{>},
      (2.3,-2);(0,-3)**\crv{ (1.5,-1.5) & (.5,-2) }?(.62)="y1",
      "y1"+(.04,.015);"y1"**\crv{}*\dir{>},
      "stop";"TopArcStop"**\crv{"stop"+(1,.4) &  (4,2) &"TopArcStop"+(1,1.5)},
    \endxy}$
    &
    $\vcenter{\xy /r9pt/:,
      (-2,0);(2,0)**\crv{ (-1,.5) & (0,2) & (-1,3) & (-2,2) & (-1,.5) & (0,0)
       & (1,-.5) & (2,-2) & (1,-3) & (0,-2) & (1,-.5) }?(.24)="x"?(.8)="x1",
       "x"+(.04,.01);"x"**\crv{}*\dir{>},
       "x1"+(.04,.02);"x1"**\crv{}*\dir{>},
      (-2.3,2);(0,3)**\crv{ (-1.5,1.5) & (-.5,2) }?(.42)="y",
      "y"+(.04,.015);"y"**\crv{}*\dir{>},
      (2.3,-2);(0,-3)**\crv{ (1.5,-1.5) & (.5,-2) }?(.62)="y1",
      "y1"+(.04,.015);"y1"**\crv{}*\dir{>},
      (-2,0);(0,-3)**\crv{(-2,0)+(-1,-.4) &  (-4,-2) & (0,-3)+(-1,-1.5)},
    \endxy}$
  \\
   \noalign{\vskip4pt}
  (i) & (ii)  
  \end{tabular}
  \caption{}
\label{2 drots b}
\end{figure}

\noindent At this stage, it is known that the knot does not have any
tight $drots$ 2-groups, and so this could not occur if only one of the
orbiters was flyped. We therefore consider the possibilities that arise when
both are flyped.

In Figure \ref{2 drots c}, one min-tangle, $T$, of the orbit of orbiter
loner $c$ is shown (by means of a dotted outline). The orbiter loner $d$
is in the min-tangle $T$ of $c$'s orbit, so $d$ can only have flype positions
within $T$. Thus if $c$ and $d$ are able to flype so as to participate in an
$OTS$ 6-tangle, the only possible flype position for $c$ would be the pair of
arcs $e_5$ and $e_6$ on the side of $T$ opposite to that of $c$, and then $d$
must have the opportunity to flype to a position on arcs $e_2$ and $e_3$ where
involves exactly one of arcs $e_2$ and $e_3$ coincides with one of the arcs
$e_5$ and $e_6$ (it is not possible for $d$ to
have a flype position on the pair of arcs $e_5$ and $e_6$, since if it did,
then $c$ and $d$ would belong to the same group). In Figure \ref{2 drots c},
that part of the orbit of $d$ which consists of all min-tangles between $d$ and
the position on arcs $e_2$ and $e_3$ is indicated by the tangle (also outlined
by a dotted curve) which lies immediately above $d$. There are four cases in
total to consider, $e_2$ coincides with either $e_5$ or $e_6$, or else $e_3$
coincides with either $e_5$ or $e_6$.

\begin{floatingfigure}{.3\hsize}
 \hbox {\hfil$\vcenter{\xy /r11pt/:,
      (-2,0);(1.5,-.5)**\crv{ (-1,.5) & (0,2) & (-1,3) & (-2,2) & (-1,.5) & (0,0)
       & (1,-.3) & (2,-2) & (1,-3) & (0,-2) & (1,-.5) }?(.24)="x"?(.8)="x1",
       "x"+(.04,.01);"x"**\crv{}*\dir{>},
       "x1"+(.04,.02);"x1"**\crv{}*\dir{>},
      (-2.3,2);(0,3)**\crv{ (-1.5,1.5) & (-.5,2) }?(.42)="y",
      "y"+(.04,.015);"y"**\crv{}*\dir{>},
      (2.3,-2);(0,-3)**\crv{ (1.5,-1.5) & (.5,-2) }?(.62)="y1",
      "y1"+(.04,.015);"y1"**\crv{}*\dir{>},
      (1.9,-.75)*{\hbox{\small $c$}},
      (-1,0)*{\hbox{\small $d$}},
      (1.5,-.5);(1,3)**\crv{ (2,-.3) & (2.5,1) & (1,2) },
      (-2,0);(-3.2,3)**\crv{ (-2.5,-.3) & (-3.5,1.5) & (-3.2,2.5) },
      (-1,1.25);(-1,4.25)**\crv{~*=<2.5pt>{.} (-.1,1.25) & (.3,2) & (.6,2.75)
           & (.3,3.5) & (-.1,4.25)}?(.6),
      (-1,1.25);(-1,4.25)**\crv{~*=<2.5pt>{.} (-2.2,1.25) & (-2.5,2) & (-2.8,2.75)
           & (-2.55,3.5) & (-2.1,4.25)},
      (-1,-1);(-1,6)**\crv{~*=<2.5pt>{.} (0,-1) & (2,1) & (2.1,2.5) & (2,4) & (0,6) }?(.6)="t",
      (-1,-1);(-1,6)**\crv{~*=<2.5pt>{.} (-2,-1) & (-4,1) & (-4.1,2.5) & (-4,4) & (-2,6) },
      (1.2,-1.25);(1.2,-4.25)**\crv{~*=<2.5pt>{.} (2.2,-1.25) & (2.6,-2) & (2.9,-2.75)
           & (2.6,-3.5) & (2.2,-4.25)},
      (1.2,-1.25);(1.2,-4.25)**\crv{~*=<2.5pt>{.} (.2,-1.25) & (-.2,-2) & (-.5,-2.75)
           & (-.2,-3.5) & (.2,-4.25)},
      (-.6,5.6);(1.7,-3.8)**\crv{ (-.6,7.5) & (1.5,7) & (3.9,2) & (4,-1) & (3.2,-5) & (1.7,-5.9)},     
      (-1.6,5.6);(.7,-3.8)**\crv{ (-1.6,7.5) & (-3.7,7) & (-5.4,2) & (-4.5,-1)
      & (-2.8,-5) & (.7,-5.9)},
      (-1.6,3.9);(-1.6,4.5)**\dir{-},
      (-.6,3.9);(-.6,4.5)**\dir{-},
      (-1.6,4.5)*!<5pt,0pt>{\hbox{\small $e_2$}},
      (-.6,4.5)*!<-6pt,0pt>{\hbox{\small $e_3$}},
      (-3.2,3)*!<4pt,0pt>{\hbox{\small $e_1$}},      
      (1,3)*!<-5pt,0pt>{\hbox{\small $e_4$}},      
      (-1.6,6.1)*!<5pt,-1pt>{\hbox{\small $e_5$}},
      (-.6,6.1)*!<-6pt,-1pt>{\hbox{\small $e_6$}},
      "t"*!<-6pt,0pt>{\hbox{\small $T$}},
   \endxy}$\hfil}
  \caption{}
\label{2 drots c}
\end{floatingfigure}

Consider the possibility that $e_2$
coincides with $e_5$. Then $e_3$ can't coincide with $e_6$, and $e_1$ can't
coincide with (that is, be connected directly to) $e_6$, since if it were, the
$e_1$-$e_6$ arc would form a crossing with the $e_2$-$e_5$ arc, which is not
possible since $T$ is a min-tangle, and the two arcs on one side of a min-tangle
do not meet in a crossing inside the min-tangle. Thus $e_1$ must cross the
$e_2$-$e_5$ arc to meet up with arc $e_3$ and therefore arc $e_4$ must coincide
with arc $e_6$. But this implies the existence of a 2-tangle inside $T$, which is
not possible. Next, consider the possibility that $e_2$ coincides with $e_6$.
Then $e_1$ does not coincide with $e_5$, nor can $e_4$ coincide with $e_5$,
otherwise arcs $e_5$ and $e_6$ form a crossing inside $T$ and this is not possible.
Thus $e_5$ must coincide with $e_1$ and so $e_3$ must coincide with $e_4$ (note that
in this case, $e_3$ cannot be the arc crossing the loop to form the positive $drots$
2-group for which $d$ is the loner orbiter, since that would imply that the knot had
a tight $drots$ 2-group). But then
arc $e_1$-$e_5$ will not meet arc $e_3$-$e_4$ after $d$ is flyped to lie on arcs
$e_2$ and $e_3$, and thus when $c$ is flyped to lie on arcs $e_5$ and $e_6$,
$c$ and $d$ will be joined by arc $e_2$-$e_6$, with arcs $e_1$-$e_5$ and  $e_3$-$e_4$
being the other two arcs that would have to belong to an $OTS$ 6-tangle involving
$c$ and $d$ in their new positions. Thus this case also does not lead to any
$OTS$ 6-tangles involving $c$ and $d$. The remaining two cases are similar, with
none resulting in an $OTS$ 6-tangle involving $c$ and $d$.

Thus any of the master arrays that we are still considering that contain exactly 
two positive $drots$ 2-groups with a configuration as shown in Figure \ref{2 drots
 a} (i) will not need to be submitted to $OTS$.

\begin{floatingfigure}{.2\hsize}
 \hbox to .2\hsize{\hfil$\vcenter{\xy /r40pt/:,
  (.75,0)="a1",
  (1,0)="a2",
  (.125,1);(.625,1.375)**\crv{ (.275,.9) & (.475,1.1) }?(.42)="y",
     "y"+(.04,.035);"y"**\crv{}*\dir{>},
  (1.125,1.375);(1.625,1)**\crv{ (1.275,1.1) & (1.475,.9) }?(.32)="y",
     "y"+(.04,-.04);"y"**\crv{}*\dir{>},      
  "a1";"a2"**\crv{(.75,.25) & (1,.5) & (1,1) & (1.375,1.375) & (1.625,1.375) & (1.625,1.125) & (1.25,.75) & (.5,.75) & (.125,1.125) & (.125,1.375) & (.375,1.375) & (.75,1) & (.75,.5) & (1,.25)}?(.7)*\dir{}="c"+<16pt,.75pt>;"c"**\dir{}*\dir{>}, 
  (.875,.4)*!<-6pt,0pt>{\hbox{\small $b$}}, 
  (.75,.8)*!<6pt,5pt>{\hbox{\small $c$}}, 
  (1,.8)*!<-6pt,5pt>{\hbox{\small $d$}}, 
  \endxy}$\hfil}
  \caption{}
  \label{2 drots d}
\end{floatingfigure}

We now consider a knot with exactly two positive $drots$ 2-groups for which there
is a configuration as shown in Figure \ref{2 drots a} (ii). We need only
consider those $OTS$ 6-tangles which contain each of the two orbiter loners in some
position. In particular, if the two loners are in the position as shown in 
Figure \ref{2 drots a} (ii), the only possible $OTS$ 6-tangle they could 
participate in would require a crossing as shown in Figure \ref{2 drots d},
where the crossings that make up the $OTS$ 6-tangle are labelled $b$, $c$
and $d$, with $c$ and $d$ being the two orbiter loners. Our immediate task
is to determine what other, if any, $OTS$ 6-tangles might involve $c$ or
$d$, in whatever positions they occupy of their respective orbits.
Now, both $c$ and $d$ belong to a min-tangle $T$ in the orbit of $b$. Let 
$e$ and $f$ denote the incident arcs on the side of $T$ that is opposite 
to $b$. Since neither $c$ nor $d$ can flype to arcs outside of $T$, other
than possibly  the  arcs that are incident to $T$, it suffices to examine
$T$.

 Suppose that it is possible to flype both $c$ and $d$ in such a way
that in their new positions, they have an arc in common (so that, in
their new positions, $c$ and $d$ are adjacent). Then there must exist
tangles $T_1$ and $T_2$, each consisting of one or more min-tangles from
the orbits of $c$ and $d$, respectively, with $c$ adjacent to $T_1$
and $d$ adjacent to $T_2$, and one of the incident arcs on the other side
of $T_1$ is coincident with one of the incident arcs on the other side of
$T_2$, as shown in Figure \ref{2 drots e}. 

\begin{floatingfigure}{.4\hsize}
 $\vcenter{\xy /r40pt/:,
  (.75,0)="a1",
  (1,0)="a2",
  (.125,1);(.625,1.375)**\crv{ (.275,.9) & (.475,1.1) }?(.42)="y",
     "y"+(.04,.035);"y"**\crv{}*\dir{>},
  (1.125,1.375);(1.625,1)**\crv{ (1.275,1.1) & (1.475,.9) }?(.32)="y",
     "y"+(.04,-.04);"y"**\crv{}*\dir{>},      
  "a1";"a2"**\crv{(.75,.25) & (1,.5) & (1,1) & (1.375,1.375) 
   & (1.625,1.375) & (1.625,1.125) & (1.25,.75) 
   & (.5,.75) & (.125,1.125) & (.125,1.375) 
   & (.375,1.375) & (.75,1) & (.75,.5) 
   & (1,.25)}?(.7)*\dir{}="c"+<16pt,.75pt>;"c"**\dir{}*\dir{>}, 
  (.875,.4)*!<-6pt,0pt>{\hbox{\small $b$}}, 
  (.75,.8)*!<6pt,5pt>{\hbox{\small $c$}}, 
  (1,.8)*!<-6pt,5pt>{\hbox{\small $d$}}, 
  (.875,.5);(.875,2.5)**\crv{~*=<2.75pt>{.} (1.7,.5) & (2.5,1.5) & (1.7,2.5) }?(.6)="t",
  (.875,.5);(.875,2.5)**\crv{~*=<2.75pt>{.} (.05,.5) & (-.75,1.5) & (.05,2.5) },
  "t"*!<-7pt,0pt>{\hbox{\small $T$}},
  (.33,.75);(.33,1.7)**\crv{~*=<2.75pt>{.}  (.6,.75) & (.825,1.25)  & (.6,1.7) },
  (.33,.75);(.33,1.7)**\crv{~*=<2.75pt>{.}  (.06,.75) & (-.165,1.25)  & (.06,1.7) },
  (1.43,.75);(1.43,1.7)**\crv{~*=<2.75pt>{.}  (1.7,.75) & (1.925,1.25)  & (1.7,1.7) },
  (1.43,.75);(1.43,1.7)**\crv{~*=<2.75pt>{.}  (1.16,.75) & (.925,1.25)  & (1.16,1.7) },
  (.39,1.6);(1.4,1.6)**\crv{ (.39,1.9) & (.895,2) & (1.4,1.9) },
  (.25,1.6);(.31,2)="e1"**\crv{ (.25,1.8) & (.31,1.8) },
  (1.54,1.6);(1.48,2)="f1"**\crv{ (1.54,1.8) & (1.48,1.8) },
  (.65,2.2)="e";(.65,2.8)**\dir{-},
  (1,2.2)="f";(1,2.8)**\dir{-},
  "e"*!<5pt,-3pt>{\hbox{\small $e$}},
  "f"*!<-5pt,-3pt>{\hbox{\small $f$}},
  "e1"*!<5pt,3pt>{\hbox{\small $e_1$}},
  "f1"*!<-6pt,3pt>{\hbox{\small $f_1$}},
  \endxy}$
  \caption{}
  \label{2 drots e}
\end{floatingfigure}

 Since $b$ could then be flyped to arcs $e_1$ and $f_1$,
it follows that $e_1$ must coincide with $e$ and $f_1$
must coincide with $f$, whence $c$ and $d$ in their new positions would
form an $OTS$ 6-tangle with $b$ flyped to arcs $e$ and $f$. There is
no other possible $OTS$ 6-tangle that $c$ and $d$ (in any of their
positions) can participate in. For such a master array, it suffices to 
consider only $OTS$ 6-tangles which involve $c$, $d$ and the crossing $b$, 
and there are at most two such $OTS$ 6-tangles.

There is a special case of this last scenario. If the arc that cuts the
loop to form the positive $drots$ 2-group of $d$ proceeds directly to
cut the loop to form the positive $drots$ 2-group of $c$, as shown in
Figure \ref{2 drots f} (i), then it is not necessary to send the master array
to $OTS$. 

 For the result of performing $OTS$ on the $OTS$ 6-tangle
which involves the two orbiter loners is shown in Figure \ref{2 drots f}
(ii). However, the tangle shown in Figure \ref{2 drots f} (v) contains
a negative 2-group, so it will already have been produced by $D$. The
results of applying three $OTS$ operations in succession (and the 
reductions that we have implemented above for $OTS$ will not prevent 
these $OTS$ operations from being performed) to the tangle shown in
(v) are displayed in (iv), (iii) and (ii), respectively, so the $OTS$
operation on the $OTS$ 6-tangle shown in (i) does not need to be
performed.

\begin{figure}[ht]
\centering
  \begin{tabular}{c@{\hskip20pt}c@{\hskip20pt}c@{\hskip20pt}c@{\hskip20pt}c}   
   $\vcenter{\xy /r25pt/:,
    (.75,0)="a1",
    (1,0)="a2",
    (.125,.875)="a3",
    (1.625,.875)="a6",
    "a3";"a6"**\crv{(.625,1.25) & (1.125,1.25)}
        ?(.15)*\dir{}="c"+<8pt,3pt>;"c"**\dir{}*\dir{>},
    "a1";"a2"**\crv{(.75,.25) & (1,.5) & (1,1) 
      & (1.375,1.375) & (1.625,1.375) & (1.625,1.125) 
      & (1.25,.75) & (.5,.75) & (.125,1.125) & (.125,1.375) 
      & (.375,1.375) & (.75,1) & (.75,.5) 
      & (1,.25)}?(.7)*\dir{}="c"+<16pt,.75pt>;"c"**\dir{}*\dir{>}, 
    \endxy}$
    &
    $\vcenter{\xy /r25pt/:,
     (.875,0)="a1",
     (1.375,0)="a2",
     (.125,1.25)="a3",
     (2.125,1.25)="a6",
     "a3";"a6"**\dir{-}?(.45)*\dir{<},
     "a1";"a2"**\crv{(.875,.5) & (1.5,1.5) & (1.75,1.5) 
        & (1.875,1.25) & (1.75,.5) & (1.5,.125) 
        & (.75,.125) & (.5,.5) & (.375,1.25) & (.5,1.5) 
        & (.75,1.5) & (1.375,.5)}
        ?(.5)*\dir{}="c"-<10pt,-.5pt>;"c"**\dir{}*\dir{>}, 
     \endxy}$
     &
     $\vcenter{\xy /r25pt/:,
       (.25,0)="a1",
       (1.75,0)="a2",
       (0,1)="a3",
       (2,1)="a6",
       "a3";"a6"**\dir{-}?(.5)*\dir{>},
       "a1";"a2"**\crv{(.25,.25) & (.75,1.25) 
           & (1,1.5) & (1.5,1.5) & (1.75,1.25)
           & (1.75,.75) & (1.25,.25) & (.75,.25) 
           & (.25,.75) & (.25,1.25) & (.5,1.5) & (1,1.5) 
           & (1.25,1.25) & (1.75,.25)}
            ?(.5)*\dir{}="c"-<10pt,-2pt>;"c"**\dir{}*\dir{>}, 
       \endxy}$
       &
       $\vcenter{\xy /r25pt/:,
         (.25,0)="a1",
         (1.75,0)="a2",
         (0,.25)="a3",
         (2,1)="a6",
         "a3";"a6"**\crv{(.5,.25) & (1,1)}?(.5)*\dir{}="c"+<10pt,6pt>;"c"**\dir{}*\dir{>},
         "a1";"a2"**\crv{(.25,.25) & (.75,1.25) 
            & (1,1.5) & (1.5,1.5) & (1.75,1.25) 
            & (1.75,.75) & (1.25,.25) & (.75,.25) 
            & (.25,.75) & (.25,1.25) & (.5,1.5) 
            & (1,1.5) & (1.25,1.25) & (1.75,.25)}
            ?(.5)*\dir{}="c"-<10pt,-2pt>;"c"**\dir{}*\dir{>},
         \endxy}$
         &
        $\vcenter{\xy /r25pt/:,
         (.25,0)="a1",
         (1.75,0)="a2",
         (0,.25)="a3",
         (2,.25)="a6",
         "a3";"a6"**\crv{(.5,.25) & (.75,.5) 
               & (1.25,.5) & (1.5,.25)}?(.5)*\dir{}="c"+<10pt,-2pt>;
               "c"**\dir{}*\dir{>},
         "a1";"a2"**\crv{(.25,.25) & (.75,1.25) & (1,1.5) 
            & (1.5,1.5) & (1.75,1.25) & (1.75,.75) 
            & (1.25,.25) & (.75,.25) & (.25,.75) 
            & (.25,1.25) & (.5,1.5) & (1,1.5) & (1.25,1.25) 
            & (1.75,.25)}?(.5)*\dir{}="c"-<10pt,-2pt>;"c"**\dir{}*\dir{>},
         \endxy}$ \\[10pt]
       (i) & (ii) & (iii) & (iv) & (v) 
      \end{tabular}
     \caption{}
     \label{2 drots f}
\end{figure}

In summary, if a master array contains exactly two positive $drots$
2-groups, then there is only one scenario that will cause us to submit
the master array to $OTS$, and that is the one illustrated in Figure
\ref{2 drots d}, but not as shown in Figure \ref{2 drots f} (i).
         
\noindent {\bf Configuration dependent reductions for $OTS$.}
The observation to be made here for a given master array $M$ is that
if $M$ contains a negative 2-group, then we need only consider $OTS$
6-tangles that contain a crossing from any negative 2-group of $M$,
and if $M$ contains a positive 3-group, then we need only consider $OTS$
6-tangles that contain a crossing from any positive 3-groups of $M$,
and if $M$ contains a positive $drots$ 2-group, then we need only consider
$OTS$ 6-tangles that contain the orbiter loner of any positive $drots$
2-group of $M$. Moreover, if the master array contains exactly
two positive $drots$ 2-groups, then there are only two $OTS$ 6-tangles
that we shall apply $OTS$ to, and these are the two positions of the
$OTS$ 6-tangle comprised of the loner orbiters of the two positive
$drots$ 2-groups, together with the third crossing as shown in
Figure \ref{2 drots d}. Of course, if the master array
is without any one of the above-mentioned objects, the corresponding 
requirement is not relevant. In particular, in a master array all of
whose groups are loners or positive 2-groups, then every $OTS$ 6-tangle
must have $OTS$ applied to it.

%% file: implem3.tex
\section{An Implementation of the Construction, Part III: the Operators
  at Work}

In this, the final section, we describe an efficient implementation of
the algorithms described in the preceding sections. As the construction of
the knots at a given crossing size proceeds, it is important
to keep in memory a database of the knots of that crossing size that have
been constructed in order that as each successive knot is produced, it 
can be compared against the database to determine whether
it is a new knot, in which case it is added to the database. However,
the size of such a database becomes an issue even at 18 crossings. Our first
reduction in memory requirements stems from the observation that the master
array is completely determined by any configuration of the knot, and thus
to compare master arrays, it suffices to compare the special configuration
obtained from each master array by placing each group in its zero position.
In fact, this might even be considered to be the so-called {\it ideal
configuration} of a prime alternating knot. For our internal comparisons, we
shall store the zero position configuration, and in fact, we store this in
its Dowker-Thistlethwaite code. 

 The other observation that leads to reduced database memory requirements
 derives from the fact (to be established in this section) that
 we can have the operators act on master arrays to produce master group 
 codes for the knots that result from the operator application, and there
 is a strong connection between the master array of the knot being operated
 on and the resulting master array. We shall exploit this connection to
 obtain a method for structuring the database so that only part of it
 needs to be available at any one stage of the construction process. This
 partitioning of the database is done according to the number of groups in
 a group code for the knot.

\begin{definition}\label{def of gn}
 Let $K$ be a prime alternating knot. Then the number of groups in any
 full group configuration of $K$ is called the {\em group number} of $K$,
 denoted by $\gn{K}$.
\end{definition}

 We shall begin the process of generating the prime alternating knots of
 $n+1$ crossings with the application of $DROTS$, which takes as its input
 the collection of all prime alternating knots of $n$ crossings, partitioned
 by group number. The knots that are produced by $DROTS$ are the $K_A$ knots
 of $n+1$ crossings, and each new $K_A$ knot will be assessed according to
 the reduction rules for $T$ and $OTS$ to see whether or not it should be
 put into either or both of the input queues for $T$ and for $OTS$, as well
 as being written to the $n+1$ crossing $K_A$ directory, in files according
 to group number.

 Our approach will be to set $k$ equal to the value of the smallest group
 number for the knots at $n$ crossings, and for each $n$ crossing knot
 $K$ of group number $k$, its master array $M$ is sent simultaneously to
 the four $DROTS$ operator situation; namely $D$ on negative groups and
 loners, $D$ on positive 2-groups, $ROTS$ on negative 2-groups and $ROTS$
 on negative 3-groups. Each of these operators will examine $M$ and decide,
 in view of their respective restrictions, whether or not the operator
 should be applied to one or more of the knot configurations of $K$. Any
 knot that does result from an application of one of these four operators
 is then examined to see whether or not it has already been produced. If
 not, then it is added to the database in memory, examined to see whether
 or not it should be subsequently submitted to one or both of the operators
 $T$ and $OTS$, and ultimately saved to disk, classified by its group number
 within the collection of $K_A$ knots. When the last $n$ crossing knot of
 group number $k$ has been processed, we increment $k$ and repeat the
 process.

 The effect of each of the various operators in the $DROTS$ family has on 
 the group number of a knot is easily described. To begin with, $D$ acting 
 on negative groups and loners does not change the group number of the knot
 on which it operates, while $D$ acting on positive 2-groups and $ROTS$ 
 acting on negative 2-groups increases the group number by one (the 2-group
 being operated on is converted into a loner and a 2-group, the 2-group 
 becoming a new min-tangle in the orbit of the loner). The last operator is
 $ROTS$ acting on negative 3-groups, and it results in a negative 3-group
 being replaced by a loner in the flype orbit of the original 3-group, and
 a new min-tangle added to this orbit, the min-tangle being a positive
 $drots$-tangle. Thus $ROTS$ acting on negative 3-groups causes an increase
 of two in the group number of the knot.

 The significance of the preceding observations about the change in group
 number when a knot is acted on by a $DROTS$ operator is revealed when we
 consider the $K_A$ knot database that is being kept in memory. When the
 last $n$ crossing knot of group number $k$ has been processed, we may
 purge the $K_A$ knots with group number $k$ from our database. When
 $DROTS$ has finished processing the last $n$ crossing knot, the database
 can be completely purged, for our reduction scheme assures us that no
 subsequent operation will ever produce a $K_A$ knot. For the range of
 crossing sizes that we have computed, this results in a reduction of
 approximately $80\%$ in main memory requirements.

 We have one last general remark to make. As we have structured our
 operators, once a master array has been made, it is submitted to the
 database routine. It is the responsibility of the database routine
 to determine first whether the knot is new and therefore to be added
 to the database or that it already exists and the submitted master
 array is to be discarded. If the knot is new, the database routine
 will compute the Dowker-Thistlethwaite code for the zero position
 configuration and add it to the database. As well, it will determine
 the group number of the submitted knot and write it to the appropriate
 disk file, according to type ($K_A$, $K_B$, or $K_C$) and group number.
 Furthermore, during the construction of the $K_A$ knots, the database
 routine will be responsible for building the input files for the
 applicaton of $T$ and $OTS$.
 
 At the completion of the construction, we will have all $K_A$ knots,
 all $K_B$ knots, and all $K_C$ knots, each collection classified according
 to group number.
 
\subsection{$D$ on negative groups and loners}

 $D$ receives the master array $M$ of a knot of $n$ crossings, and searches
 for negative groups and loners. For each negative group or loner with a
 flype orbit in the negative direction, a master group code
 for the knot that results from the application of $D$ to the group or loner
 in question in position zero is constructed from $M$ as follows. Let $G$
 denote the group in question. The first occurrence of a group arc of $G$ in
 position zero is located, and then the second of $G$'s arcs in position
 zero is located. The segment of the master array that lies between these
 two positions referred to as the first segment for $G$ in the master array.
 Reverse the order of the entries in the first segment, and then replace
 the label $G$ in each of its positions by the next available label for a
 group of size one greater than the size of $G$, recording it as a negative
 group. Finally, a sign change must be recorded for each of the groups in
 the core tangle whose arcs were what we called the starter core group arcs
 (those groups which have one arc in the first section and one arc in the
 second section of $G$'s orbit).

 On the other hand, if $G$ is a loner with an orbit in the positive
 direction, then a master group code for the knot that results from an
 application of $D$ to $G$ is obtained from $M$ as follows: first, delete
 all positions of $G$ except for the one in zero position, then reverse
 the first segment for $G$ in the master array and replace the label $G$
 by the next available label for a group of size two, recording it as a
 negative group, and finally, change the sign of each group that had one
 arc in the first section of the loner and one arc in the second section
 of the loner.

 In either case, the result is a master group code for the knot that
 results from the application of $D$ to $G$. The master group code can
 now be checked according to the reduction criteria described in the
 preceding section to see if we should proceed further with this group.
 If so, and $G$ is a negative group or a loner with an orbit in the
 negative direction, then we examine the orbit of $G$ for flype-symmetric
 tangles in order to identify the actual positions in the orbit of $G$ at
 which we shall apply $D$, while if $G$ is a loner with an orbit in the
 positive direction, we shall apply $D$ to $G$ in each of its positions,
 in the manner described above for $G$ in position zero. When the positions
 have been identified, we then construct from $M$ a master group code for
 each knot that is produced by applying $D$ to $G$ in each of the identified
 positions. The resulting master group codes are converted to master arrays
 and submitted to the database.

\subsection{\bf $D$ on positive 2-groups}

 Suppose that $G$ is a positive 2-group in the knot $K$ of $n$ crossings.
 We may form a master group code for the knot that results from an
 application of $D$ to $G$ in position $j$ of its orbit as follows. Let
 $m$ denote the number of loners in $K$, and suppose that $G$ was the
 $i^{th}$ 2-group to have been labelled. Further, suppose that the orbit of
 $G$ had $k$ positions. In $M$, locate both arcs of $G$ in position $j$.
 These two positions partition the orbit of $G$ into two segments, and we
 choose one and reverse it. Then replace the first label $2_{i}^j$ by
 $1_{m+1}^{j},-2_i^0,1_{m+1}^{j}$ and the second label $2_{i}^j$ by
 $1_{m+1}^{k+1},-2_i^0,1_{m+1}^{k+1}$. Next, for every $0<t\ne j$, replace
 each occurrence of $2_i^t$ by $1_{m+1}^t$. Finally, each group that
 had one arc in the first section for $G$ in the master array and one arc
 in the second section must have a change of sign recorded.

 To begin with, we apply the above process to $G$ in position zero, and
 then apply the reduction criteria for $D$ on positive 2-groups. If $G$
 the reduction criteria do not stop $D$ from being applied to $G$, then
 we do apply $D$ to $G$ in each position and submit the resulting negative
 2-group to a group $LNR$ competition against the other negative 2-groups
 (in this setting, all the negative 2-groups in the resulting knot will
 be 2-groups in a negative $drots$ tangle, and so will have trivial flype
 orbit). If the group that results from the application of $D$ to $G$ in
 a given position wins its group $LNR$ competition, then the master array
 of the knot that resulted from the application of $D$ to $G$ in the
 given position is calculated and submitted to the knot database.

\subsection{\bf $ROTS$ on negative 2 and 3-groups}

 Given the master array $M$ of an $n$ crossing knot $K$, we check first
 to see if any reduction rules apply. If $M$ contains two or more negative
 groups, or just one negative group, but that negative group either has
 size greater than 3 or else is the negative 2-group of a negative $drots$
 tangle, then $M$ is not submitted to $ROTS$.

 Once it has been determined that $M$ is to be submitted to $ROTS$, the
 location of the unique negative group $G$ in $M$ is known, and $G$ is
 either a 2-group or a 3-group. If $G$ is a 2-group, then we apply $ROTS$
 to the group once for each position of the group in its orbit. A master
 group code for the knot that results from the application of $ROTS$ to
 $G$ in position $i$ is obtained from $M$ as follows. First, delete all
 copies of $G$ in positions other than $i$. The result of applying $ROTS$
 to $G$ is to replace the $G$ by a positive 2-group and a loner. Let $m$
 denote the number of loners in $M$, and let $G=-2_k$. The first occurrence
 of $-2_k^i$ in $M$ is replaced by $1_{m+1}^0,2_k^0,1_{m+1}^0$, and the
 second occurrence of $-2_k^i$ in $M$ is replaced by
 $1_{m+1}^1,2_k^0,1_{m+1}^1$.

 It remains to deal with the case when $G$ is a 3-group. The application of
 $ROTS$ to a negative 3-group replaces the 3-group by a positive 2-group
 with trivial flype orbit, a loner with a two position flype orbit (one
 min-tangle in the orbit of this loner is the positive 2-group that was
 created), and a loner that inherits the orbit of $G$. We apply $ROTS$
 to $G$ in each position of its orbit. A master group code for the knot
 that results from the application of $ROTS$ to $G$ in position $i$ is
 obtained from $M$ as follows. Let $m$ and $t$ denote the number of loners,
 respectively 2-groups, in $M$, and suppose that $G=-3_k$. In $M$, replace
 the first occurrence of $-3_k^i$ by $1_{m+2}^0,1_{m+1}^0,2_{t+1}^0,
 1_{m+1}^0,1_{m+2}^1$, and replace the second occurrence of $-3_k^i$ by
 $1_{m+2}^1,1_{m+1}^1,2_{t+1}^0,1_{m+1}^1,1_{m+2}^0$. Replace every
 occurence of $-3_k^j$ for $j\ne i$ by $1_{m+2}^j$.

 In either case, the result is a master group code for the knot that results
 from the application of $ROTS$ to the negative group in one position of
 its orbit. The master array is then constructed from the master group code,
 and the knot is submitted to the database. 

 Once $D$ and $ROTS$ are finished, all $K_A$ knots have been constructed.
 The database can now be purged completely in readiness for the construction
 of the $K_B$ and $K_C$ knots.

 Recall that as it processed the $K_A$ knots, the database routine was
 preparing the input files for $T$ and $OTS$. In the case of $T$, each
 $K_A$ master array $M$ is checked to see if it contains a negative
 $DROTS$ 2-group, or no positive 2-groups, or two or more positive
 $DROTS$ 2-groups. If any of these conditions are met, then $M$ is not
 submitted to $T$. Thus the input file for $T$ will consist of all master
 arrays from the $DROTS$ output that contain no negative $DROTS$ 2-groups,
 at least one positive 2-group, and at most one positive $DROTS$ 2-group.

\subsection{\bf $T$ on positive 2-groups}

 Consider the master array $M$ of an $n+1$ crossing $K_A$ knot that has
 been submitted to $T$.

\noindent Case 1: $M$ contains a single positive $DROTS$ 2-group $G$. Then 
 that is the only 2-group that $T$ need be applied to. The positive $DROTS$ 
 2-group $G$ necessarily has trivial orbit, and turning $G$ will increase 
 the size of the orbiter group $H$ by 2 crossings. Furthermore, this
 enlarged group will retain the orbit of $H$. The group number of the new
 knot will be one smaller than the group number of $K$. We may readily
 construct from $M$ a master group code for the knot that results from the
 application of $T$ to $G$. To begin with, $G$ is denoted by $2_i$ for some
 index $i$, and in $M$, we shall find $2_i$ surrounded by the arcs of some
 position $r$ of its orbiter group $H$, which must be a negative group, say
 $H=-m_j$. Thus $M$ will have the form
 $$
  \ldots,-m_j^r,2_i^0,-m_j^r,\,\ldots\hbox to 0pt{\hss\xy /r20pt/:,
     (0,0)*{\null},(-.9,-.2);(-.9,-.4)**\dir{-};(.9,-.4)**\dir{-};
    (.9,-.2)**\dir{-},
     (0,-.7)*{\text{\footnotesize Section $1$}},\endxy\hss}\ldots,-m_j^s,2_i^0,-m_j^s,\ldots
 $$
 Reverse the subsequence labelled section $1$ above. Then, if there are already $q$
 groups of size $m+2$ in $K$, replace the subsequence $-m_j^r,2_i^0,-m_j^r$
 by $(m+2)_{q+1}^r$, replace the subsequence $-m_j^s,2_i^0,-m_j^s$ by
 $(m+2)_{q+1}^s$, and replace every other position $-m_j^t$ by
 $(m+2)_{q+1}^t$. Finally, for each of the remaining groups, any that had
 one arc in section 1 and one arc in the other section must have their sign
 changed. This completes the construction of a master group code for the
 knot that has been produced from $K$ by turning $G$.

 This master group code is then submitted to the procedure that will
 apply the reduction rules to determine whether to construct its master
 array and submit it to the database procedure or to discard it.

\noindent Case 2: $M$ does not contain any positive $DROTS$ 2-groups. For 
 each position of each positive 2-group $G$, we perform the following
 procedure. Suppose that $G$ is the 2-group $2_i$ and we are considering
 position $r$. In $M$, we locate all positions $2_i^s$ for $s\ne r$ and
 remove them. Then locate the subsequences contained between the first
 and second occurrences of the arcs $2_i^r$ and reverse it. Next, examine
 the code to locate any groups that exactly one of its group arcs in the
 subsequence that was reversed and change the sign of the group. Finally,
 replace the two arcs $2_i^r$ by $2_i^0$. The result is almost a master
 group code for the knot that results from turning $G$ in position $r$.
 All of the orbit information for every group except the turned group is
 already in place; that is, the array that has been constructed is that
 array that would have been built if we had started with a group code for
 the knot and constructed a master group code from it, leaving the group
 $2_i$ for last. In effect, we are picking up the process at that point.
 We complete the construction of a master group code for the knot that
 results from this application of $T$ by identifying the orbit of $2_i$.
 It is worthwhile to note here that if the group $2_i$ had a non-trivial
 orbit, then when it is turned in any position of its orbit, the resulting
 2-group will have trivial orbit, and so no further work is required;
 that is, we will already be in possession of a master group code for the
 resulting knot. 

 This master group code is then submitted to the procedure that will
 apply the reduction rules to determine whether to construct its master
 array and submit it to the database procedure or to discard it.
 We point out that in Case 2 above, if the group $2^i$ has a non-trivial
 flype orbit, it is sufficient to apply the configuration dependent
 reduction criteria to $2^i$ in any one position of its orbit. If it is
 rejected in that position, it will be rejected in every position, and
 we can bypass the application of $T$ to group $2^i$ in any position.
 Furthermore, although we have not described them here, it is evident
 that we could set up group $LNR$ competitions for $T$ on positive 2-groups
 similar in principle to those we have described for $D$.

 The procedure that applies the reduction rules will examine the master
 group codes that are submitted to it from $T$. If the master group code
 contains any negative group, the master group code is discarded.
 If the master group code contains a positive group of size greater than
 the group that results from turning $G$ (which is a 2-group if $G$ was not
 the 2-group of a positive $DROTS$ tangle, but it could have size much
 larger than 2 if $G$ was the 2-group of a positive $DROTS$ tangle), the
 master group code is discarded. Otherwise, the master array is constructed
 from the master group code, and the knot is submitted to the database. If
 the knot is not already in the database, the database routine will insert
 it, and also write it to disk as a $K_B$ knot, classified according to its
 group number. As well, it will be submitted to the routine that is
 preparing the input files for the initial $OTS$ application.

\subsection{OTS}

 Unlike the other operators, $OTS$ does not work directly on groups, but 
 rather on 6-tangles. Consequently, $OTS$ has the potential to be more
 expensive to perform on group code. We shall therefore expend considerable
 effort towards reducing the work that $OTS$ must do.  Recall that each knot
 that was produced by $DROTS$ was examined to see if it contained a negative
 $DROTS$ tangle, and only if it did not was it forwarded to the routine that
 is preparing the input file for the application of $OTS$. As well, the output
 from the application of $T$ is sent to this routine as well (no knot that is
 produced by the application of $T$ as described above will contain a
 negative $DROTS$ tangle, so it is not necessary to examine these knots
 for negative $DROTS$ tangles).

 The job of the routine that is preparing the input file for the $OTS$
 application is to apply the reduction criteria described in the preceding
 section to each submitted knot, discarding those that do not meet the
 criteria for submission to $OTS$. In the course of the examination of a
 knot for the purpose of applying these criteria, a considerable amount of
 information about the knot is gathered. This information is pertinent to
 the identification of the $OTS$ 6-tangles to which the $OTS$ operation is
 to be applied for those knots that do get submitted to $OTS$. Thus the
 input file for the first application of $OTS$ will contain not only
 the master arrays of the knots to which $OTS$ is to be applied, but as
 well, each master array will be accompanied by additional data that will
 help to identify the $OTS$ 6-tangles to which $OTS$ is to be applied.

 As each knot that is constructed by the first application of $OTS$ is
 submitted to the database routine, it will be examined to see if it is
 new and if so, it will be inserted in the database and written to disk as
 a $K_B$ knot classified by group number.
 Furthermore, it will be added to the input file for the next application
 of $OTS$, which is to say that $OTS$ is to be applied iteratively, until
 no new knots are constructed by $OTS$. Only the very first application
 of $OTS$ will produce $K_B$ knots, as all subsequent applications of
 $OTS$ to knots that have been produced by earlier $OTS$ applications
 will only produce $K_C$ knots (actually, $K_B$ knots may be constructed,
 but they will be discarded by the database routine, since we have
 retained the $K_B$ knots in the database.
 After $OTS$ has finished processing its input file, it will
 restart with the new input file that was constructed during the just-completed
 $OTS$ run. This will continue until $OTS$ encounters an empty input file.

 We emphasize again that after the first, non-iterative application of
 $OTS$, all $K_B$ knots will have been constructed, and there are no
 reduction criteria applicable for any subsequent $OTS$ applications.

 The combined output from all of the $OTS$ rounds is also written to a
 file to serve as the input for the subsequent round of application of $T$.
 The construction will continue in this way, alternating between $T$ and
 $OTS$, with the output from each operator serving as the input for the
 next round of the other operator. The process stops when no new knots
 are produced during the application of one of the operators.

 Let us examine how the routine that will prepare the input file for
 the first, non-iterative application of $OTS$ to the knots that have been
 produced by $DROTS$ and $T$. Each master array that is submitted to the
 routine will be examined group by group. In doing so, counts are
 maintained of the number of negative 2-groups that are not negative
 $drots$ 2-groups, the number of positive 3-groups,
 and the number of positive $DROTS$ 2-groups with loner orbiter group that
 have been encountered, and for each such object, a record of that object
 and its index in the master array is constructed. If at any point in the
 examination, we encounter a negative group of size at least three, a
 positive group of size at least 4, a positive $DROTS$ 2-group whose
 orbiter group has size at least 2, or a tight $DROTS$ 2-group, or the
 number of positive $DROTS$ 2-groups with loner orbiter exceeds 2, or the
 sum of the number of negative 2-groups that are not $drots$ 2-groups, the
 number of positive 3-groups  and the number of positive $DROTS$ 2-groups
 with loner orbiter exceeds 4, then the master array is discarded.

 If we reach the end of the examination of the master array without
 discarding it, then we know that the master array array contains at most
 2 positive $DROTS$ 2-groups. For those that do contain 2 positive $DROTS$
 2-groups, we perform an additional check to see if the 2 positive
 $DROTS$ 2-groups are in the configuration shown in Figure \ref{2 drots d},
 but not as shown in Figure \ref{2 drots f} (i). If not, then the master
 array is discarded.

 For a master array $M$ that survives the application of the reduction rules,
 the next step is to identify the actual $OTS$ 6-tangles in the master
 array to which $OTS$ is to be applied. The following observations are
 pertinent. For any negative 2-group or positive 3-group $G$ of $M$, it is
 necessary only to consider $OTS$ 6-tangles that contain a crossing of $G$,
 and for any positive $DROTS$ 2-group $H$ of $M$, it is necessary only to
 consider $OTS$ 6-tangles that contain the loner orbiter of $H$. Moreover,
 if $M$ contains two positive $DROTS$ 2-groups, then we shall apply $OTS$
 to only two $OTS$ 6-tangles; namely the one that consists of the loner
 orbiters of the two positive $DROTS$ 2-groups together with the third
 crossing as shown in Figure \ref{2 drots d}, and the one obtained by
 putting the two orbiter loners in their other position, and flyping the
 third crossing to the position in its orbit that has it make an $OTS$
 6-tangle with the two loner orbiters again. If the groups in $M$ are
 either loners or positive 2-groups, and none of the positive 2-groups
 are positive $DROTS$ 2-groups, then we must apply $OTS$ to every $OTS$
 6-tangle.

 There remains the issue of identifying the $OTS$ 6-tangles in the
 master array. Of course, in the first application of $OTS$, we may
 have many restrictions on the $OTS$ 6-tangles to which $OTS$ is to be
 appliced, while in subsequent applications, every $OTS$ 6-tangle must
 have $OTS$ applied to it. The master array is a master group code, and
 an $OTS$ 6-tangle can involve as many as 3 different groups.
 It is important to be able to work with the master array rather than
 an individual group code, since any flype operation that does not involve
 any of the crossings of an $OTS$ 6-tangle will commute with $OTS$ applied
 to that $OTS$ 6-tangle. Accordingly, if we apply $OTS$ to the same $OTS$
 6-tangle in two configurations that differ only by flypes that do not involve the
 crossings of the $OTS$ 6-tangle, the resulting knots will be flype-equivalent.
 In such a case, it is not necessary to perform $OTS$ on both of the initial
 configurations. Moreover, the orbit structure of the knot that results from an
 application of $OTS$ will be very similar to the orbit structure of the
 original knot. Since we do not want to rediscover information that is
 already in our possession, we shall identify the
 various $OTS$ 6-tangles in the master array to which we are to apply $OTS$,
 and for each such 6-tangle, we construct a master group code for the knot
 that is obtained upon the application of $OTS$ to the selected $OTS$ 6-tangle.
 Just as was the case for $T$, we shall find ourselves with most of the orbit
 structure already known to us, with the orbits of at most three groups
 to be determined before we have a master group code for the knot that results
 from the application of $OTS$.

 We now describe the process whereby we are able to identify $OTS$ 6-tangles
 in the master array, and to describe the effect of applying $OTS$ to a
 selected $OTS$ 6-tangle in the master array. As we mentioned above, the
 result is a group code that is most of the way through the process of being
 converted to a master group code for the newly constructed knot. We then
 describe how to identify the groups that still need to have their orbits
 determined, which then allows us to finish the job of constructing a master
 group code for the newly constructed knot.

\begin{figure}[ht]
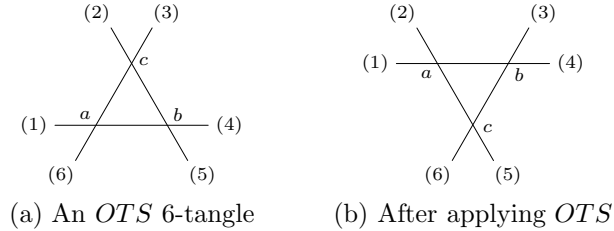

  \centering
    \begin{tabular}{c@{\hskip30pt}c}
    \vtop{\xy /r30pt/:,
        {\xypolygon12"a"{~<{}~>{}}},
        "a0";"a0"+(0,-1)**\dir{}?!{"a8";"a11"}="x1",
        "a0";"a0"+(0,-1)**\dir{}?!{"a9";"a10"}="x2",
        "x2"-"x1"="z",
        "a7"+"z"="1";"a12"+"z"="4"**\dir{-}?(.25)*!<3pt,-4pt>{\hbox{$\ssize a$}},
        "a7"+"z";"a12"+"z"**\dir{}?(.75)*!<-3pt,-4pt>{\hbox{$\ssize b$}},
        "a8"+"z"="6";"a3"+"z"="3"**\dir{-}?(.75)*!<-4pt,0pt>{\hbox{$\ssize c$}},
        "a11"+"z"="5";"a4"+"z"="2"**\dir{-},
        "1"*!<8pt,0pt>{\hbox{$\ssize (1)$}},
        "2"*!<5.5pt,-5.5pt>{\hbox{$\ssize (2)$}},
        "3"*!<-5.5pt,-5.5pt>{\hbox{$\ssize (3)$}},
        "4"*!<-8pt,0pt>{\hbox{$\ssize (4)$}},
        "5"*!<-5.5pt,5.5pt>{\hbox{$\ssize (5)$}},
        "6"*!<5.5pt,5.5pt>{\hbox{$\ssize (6)$}},        
    \endxy}
    &
    \vtop{\xy /r30pt/:,
        {\xypolygon12"a"{~<{}~>{}}},
        "a6"="1";"a1"="4"**\dir{-}?(.25)*!<3pt,4pt>{\hbox{$\ssize a$}},
        "a6";"a1"**\dir{}?(.75)*!<-3pt,4pt>{\hbox{$\ssize b$}},
        "a9"="6";"a2"="3"**\dir{-},
        "a5"="2";"a10"="5"**\dir{-}?(.75)*!<-5pt,0pt>{\hbox{$\ssize c$}},
        "1"*!<8pt,0pt>{\hbox{$\ssize (1)$}},
        "2"*!<5.5pt,-5.5pt>{\hbox{$\ssize (2)$}},
        "3"*!<-5.5pt,-5.5pt>{\hbox{$\ssize (3)$}},
        "4"*!<-8pt,0pt>{\hbox{$\ssize (4)$}},
        "5"*!<-5.5pt,5.5pt>{\hbox{$\ssize (5)$}},
        "6"*!<5.5pt,5.5pt>{\hbox{$\ssize (6)$}},        
    \endxy}\\
    \noalign{\vskip 6pt}
    (a) An $OTS$ 6-tangle & (b) After applying $OTS$
   \end{tabular}
   \caption{$OTS$ in action}
   \label{ots application}
  \end{figure}

 To begin with, we shall describe the procedure for identifying an $OTS$
 6-tangle when given a Gauss code for a knot configuration $C(K)$, following
 which, we describe the changes that must be made to this Gauss code in order
 to construct a Gauss code for the knot configuration that results from the
 application of $OTS$ to a selected $OTS$ 6-tangle in $C(K)$. To find an
 $OTS$ 6-tangle when given a Gauss code for a knot configuration, we work
 our way through the code, checking each arc to see if the two crossings that
 determine the arc participate in an $OTS$ 6-tangle. Suppose that we are
 considering the arc formed by crossings $a$ and $b$, with $b$ following
 $a$ in the Gauss code. We are looking to see if there is a crossing $c$ as
 shown in Figure \ref{ots application} (i). If so, the three crossings form
 an $OTS$ 6-tangle, and to perform $OTS$ on it, we have chosen to move the
 arc determined by $a$ and $b$ (recall that the outcomes of moving any one
 of the three arcs determined by the three crossings in the $OTS$ 6-tangle
 are all identical). The result of moving the arc $ab$ is shown in
 Figure \ref{ots application} (ii).

  We shall consider the (cyclic) Gauss code to begin with the appearance
  of the crossing $a$ that is immediately followed by crossing $b$. We scan
  the rest of the code, searching for the second occurrences of crossings
  $a$ and $b$. Once they have been found, we examine the crossings on either
  side of the second occurrence of crossing $a$ and on either side of the
  second occurrence of crossing $b$ to see if there is a crossing $c$ which
  has one of its occurrences adjacent to the second occurrence of $a$ and
  its other occurrence adjacent to the second occurrence of $b$. In the knot
  traversal, strand $(4)$ could connect next to any of strands $(2)$, $(3)$,
  $(5)$ or $(6)$, and we examine each of these possibilities. For example,
  if strand $(4)$ connects next to strand $(5)$, then there are two
  subcases to consider, depending on whether strand $(2)$ connects next to
  strand $(3)$ or to strand $(6)$. In the event that strand $(2)$ connects
  to strand $(3)$, we find that the code, which begins with $a$, $b$,
  continues with a segment of code that we shall refer to as $S_1$, then
  we find the pair $b$, $c$, followed by a segment of code that we shall
  refer to as $S_2$, then the pair $c$, $a$, followed by the final segment
  of code, which we shall refer to as $S_3$. Upon consulting Figure \ref{ots
  application}, we see that after performing the $OTS$ operation, the code
  becomes $ab$, $S_1$, $ca$, $S_2$, $bc$, $S_3$. We present the corresponding
  information in tabular form for all four cases, each with two subcases.

  \vskip4pt
    \begin{figure}[ht]
  \centering
  \begin{tabular}{r@{\hskip 15pt}lllrl}
  (a) & $(4)\to (5)$ & $(2)\to (3)$ & Before $OTS$ & $ab, S_1,bc,S_2,ca,S_3$\\
      &          &              & After $OTS$ & $ab, S_1,ca,S_2,bc,S_3$\\
  (b) &              & $(2)\to (6)$ & Before $OTS$ & $ab, S_1,bc,S_2,ac,S_3$\\
      &          &              & After $OTS$ & $ab, S_1,ca,S_2,cb,S_3$\\
  (c) & $(4)\to (6)$ & $(3)\to (2)$ & Before $OTS$ & $ab, S_1,ac,S_2,cb,S_3$\\
      &          &              & After $OTS$ & $ab, S_1,cb,S_2,ac,S_3$\\
  (d) &              & $(3)\to (5)$ & Before $OTS$ & $ab, S_1,ac,S_2,bc,S_3$\\
      &          &              & After $OTS$ & $ab, S_1,cb,S_2,ca,S_3$\\
  (e) & $(4)\to (2)$ & $(5)\to (3)$ & Before $OTS$ & $ab, S_1,cb,S_2,ca,S_3$\\
      &          &              & After $OTS$ & $ab, S_1,ac,S_2,bc,S_3$\\
  (f) &              & $(5)\to (6)$ & Before $OTS$ & $ab, S_1,cb,S_2,ac,S_3$\\
      &          &              & After $OTS$ & $ab, S_1,ac,S_2,cb,S_3$\\
  (g) & $(4)\to (3)$ & $(6)\to (2)$ & Before $OTS$ & $ab, S_1,ca,S_2,cb,S_3$\\
      &          &              & After $OTS$ & $ab, S_1,bc,S_2,ac,S_3$\\
  (h) &              & $(6)\to (5)$ & Before $OTS$ & $ab, S_1,ca,S_2,bc,S_3$\\
      &          &              & After $OTS$ & $ab, S_1,bc,S_2,ca,S_3$
  \end{tabular}
   \caption{Performing $OTS$ operations on a Gauss code}
   \label{ots scenarios}
  \end{figure}
  
 From the data in Figure \ref{ots scenarios}, we extract the following rule,
 which describes in every case how to perform the $OTS$ operation. The pair
 of crossings between
 segments $S_1$ and $S_2$ and the pair of crossings between segments $S_2$ and
 $S_3$ change places, with the crossing in common, $c$ being placed in such
 a way that it does not occupy the same relative position (first or second
 crossing) between $S_1$ and $S_2$ as it did in the original pair between
 $S_1$ and $S_2$, nor does it occupy the same relative position (first or
 second crossing) between $S_2$ and $S_3$ as it did in the original pair
 between $S_2$ and $S_3$.

 Let us now examine a group code for a knot, or more generally, the master
 array for the knot, as it provides a copy of each group in every position
 that the group can be placed. Our first objective is to identify the $OTS$
 6-tangles in the master array, and then for each $OTS$ 6-tangle, we explain
 how to transform the master array into a master group code for the knot that
 results from the application of $OTS$ to the selected $OTS$ 6-tangle. In the
 search for an $OTS$ 6-tangle, we shall examine adjacent groups $G_1$ and
 $G_2$ in the master array, say with crossing $c_1$ at the end of $G_1$
 adjacent to crossing $c_2$ at the end of $G_2$, and we look to see if there
 is a third group $G_3$ such that of the two arcs at one end of $G_3$,
 leaving crossing $c_3$ of the third group, say, one is incident to $c_1$
 and the other is incident to $c_2$. If so, then $c_1$ and $c_2$ form an
 $OTS$ 6-tangle with $c_3$, otherwise $c_1$ and $c_2$ do not participate in
 an $OTS$ 6-tangle. To perform the $OTS$ operation once we have found such
 crossings $c_1$, $c_2$ and $c_3$, we begin by creating three loners in the
 master group code: one at the end of $G_1$ next to $c_1$ (this loner is
 actually $c_1$), (this loner is actually $c_2$), and one at the end of
 $G_3$ adjacent to $c_3$ (this loner is actually $c_3$). Note that each
 crossing will actually generate two entries in the master group code,
 one for each arc of the group in question. Next, we decrement the crossing
 size of each position of $G_1$, $G_2$ and $G_3$ (using the next available
 label for groups of the respective sizes), and if this results in a group
 with no crossings (that is, we had started with a loner), then we simply
 remove all positions of that group. We then treat the crossings
 $c_1$, $c_2$ and $c_3$ just as if they were crossings in a Gauss code, and
 performs the $OTS$ operation. 
 For any two of the three crossings of the $OTS$ 6-tangle,
 if the two adjacent arcs of the $OTS$ 6-tangle that are incident to the two
 selected crossings are coincident at their other end, being arcs at the end
 of a group $G_4$, then after the $OTS$ operation, the
 third crossing of the $OTS$ 6-tangle becomes another crossing of the group
 $G_4$. Accordingly, we must relabel group $G_4$ to indicate that its size
 is now one greater than it was prior to the $OTS$ operation. The orbit
 of this enlarged copy of $G_4$ is unchanged, so $G_4$ must have its size
 increased by 1 in every position, using the next available label for a
 group of that size.

 On the other hand, if the two adjacent arcs of the $OTS$ 6-tangle that are
 incident to the two selected crossings are not coincident at their other
 ends, then the third crossing becomes a loner after the $OTS$ operation.
 Note that if any of $G_1$, $G_2$ or $G_3$ is not a loner, then the crossings
 from the other two must become loners after the $OTS$ has been performed.
 Accordingly, if two or more of $G_1$, $G_2$ and $G_3$ are not loners, then
 all three crossings $c_1$, $c_2$ and $c_3$ become loners after the
 application of $OTS$.
 
 Once we apply $OTS$, taking care of the possible enlarged groups that might
 arise, the result is almost a master group code for the knot that has been
 obtained by the application of $OTS$. All that remains is to determine the
 orbits of any (at most three) loners that were created by the application
 of $OTS$. Thus after applying $OTS$ to the master array, a very small
 amount of additional work will produce a master group code for the
 resulting knot.

 When we are examining the master array for $OTS$ 6-tangles, we shall
 take advantage of the information that we have compiled for the master
 array during the preparation of the input file for the initial round of
 $OTS$. It will often be the case that, as a result of the $OTS$ reductions,
 an $OTS$ 6-tangle will not actually have $OTS$ performed on it.
 Recall that for a given master array $M$, if $M$ contains a negative
 2-group, then we need only consider $OTS$ 6-tangles that contain a
 crossing from any negative 2-group of $M$, and if $M$ contains a
 positive 3-group, then we need only consider $OTS$ 6-tangles that
 contain a crossing from any positive 3-groups of $M$, and if $M$ contains
 a positive $DROTS$ 2-group, then we need only consider $OTS$ 6-tangles
 that contain the orbiter loner of any positive $DROTS$ 2-group of $M$.
 Moreover, if the master array contains exactly two positive $DROTS$
 2-groups, then there are only two $OTS$ 6-tangles that we shall apply
 $OTS$ to, and these are the two positions of the $OTS$ 6-tangle comprised
 of the loner orbiters of the two positive $DROTS$ 2-groups, together with
 the third crossing as shown in Figure \ref{2 drots d}.

 Of course, any knot that is produced from the initial round of $OTS$
 will consist entirely of loners or positive 2-groups, with no positive
 $drots$ tangles, so every $OTS$ 6-tangle will have to have $OTS$ applied
 to it.  

 Once the master group code has been contructed during an application of
 $OTS$, it is standardized to produce the master array, which is then
 submitted to the database routine. The database routine will compare the
 master array to the knot database, and if the knot is a new one,
 it will be inserted into the database in memory, written to
 the directory of $K_B$ or $K_C$ knots, as appropriate, organized according
 to group number, and written to the input file for the subsequent round
 of $OTS$.
 
 We continue applying the cycle of alternating iterative $OTS$, followed by
 $T$, with each operator being applied only to
 the new knots that were produced by the preceding application of the other
 operator, until no new knots are produced. At this point, the construction
 is complete.